\documentclass[11pt]{article}
\usepackage[nosumlimits]{amsmath}
\usepackage{amssymb}
\usepackage{amsthm}

\input xy
\xyoption{all}
\UseComputerModernTips

\DeclareMathOperator{\Index}{Index}
\DeclareMathOperator{\Ker}{Ker}
\DeclareMathOperator{\im}{\mathrm{Im}\,}
\DeclareMathOperator{\Coker}{Coker}
\DeclareMathOperator{\Vect}{Vect}
\DeclareMathOperator{\ch}{ch}
\DeclareMathOperator{\td}{td}
\DeclareMathOperator{\tind}{t-ind}
\DeclareMathOperator{\aind}{a-ind}
\DeclareMathOperator{\ind}{ind}
\DeclareMathOperator{\Ad}{Ad}
\DeclareMathOperator{\Aut}{Aut}
\DeclareMathOperator{\Hom}{Hom}
\DeclareMathOperator{\Gr}{Gr}

\newcommand{\tensorhat}{\mathop{\hat{\otimes}}}
\newcommand{\Cm}{ {{}^{m}C} }
\newcommand{\R}{\mathbb{R}}
\newcommand{\Q}{\mathbb{Q}}
\newcommand{\C}{\mathbb{C}}
\newcommand{\Z}{\mathbb{Z}}
\newcommand{\pt}{\mathrm{pt}}
\newcommand{\g}{\mathfrak{g}}
\newcommand{\h}{\mathfrak{h}}

\theoremstyle{plain}
\newtheorem*{theorem}{Theorem}
\newtheorem*{corollary}{Corollary}
\newtheorem*{proposition}{Proposition}
\newtheorem*{claim}{Claim}
\newtheorem*{lemma}{Lemma}

\theoremstyle{remark}
\newtheorem*{example}{Example}

\theoremstyle{definition}
\newtheorem*{definition}{Definition}
\newtheorem*{property}{Property}

\title{$K$-Theory and Elliptic Operators}
\author{Gregory D. Landweber
\bigskip\bigskip\bigskip \\ 
{Mathematics Department} \\
{University of Oregon} \\
{Eugene, OR 97403-1222} \bigskip \\
\textit{E-mail address:} \texttt{greg@math.uoregon.edu} \\
\textit{URL:} \texttt{http://math.uoregon.edu/\~{}greg} \bigskip\bigskip\bigskip
}

\begin{document}

\setcounter{section}{-1}

\maketitle

\bigskip

\begin{abstract}
This expository paper is an introductory text on topological $K$-theory and the Atiyah-Singer index theorem, suitable for graduate students or advanced undegraduates already possessing a background in algebraic topology. The bulk of the material presented here is distilled from Atiyah's classic ``$K$-Theory'' text, as well as his series of seminal papers ``The Index of Elliptic Operators'' with Singer. Additional topics include equivariant $K$-theory, the $G$-index theorem, and Bott's paper ``The Index Theorem for Homogeneous Differential Operators''. It also includes an appendix with a proof of Bott periodicity, as well as sketches of proofs for both the standard and equivariant versions of the $K$-theory Thom isomorphism theorem, in terms of the index for families of elliptic operators. A second appendix derives the Atiyah-Hirzebruch spectral sequence. This text originated as notes from a series of lectures given by the author as an undergraduate at Princeton. In its current form, the author has used it for graduate courses at the University of Oregon.
\end{abstract}

\vfill


\newpage

\tableofcontents

\newpage

\section{The Index Problem}

\subsection{Introduction}

The focus of this text is the Atiyah-Singer index theorem. One of the most significant mathematical results of the second half of this century, the index theorem provides a fundamental connection between algebraic topology, differential geometry, and analysis. In particular, it expresses the index of linear elliptic differential operators in terms of certain topological invariants.  The problem of finding such an expression was originally posed in a 1960 paper by Israel Gel'fand, and it was solved by Michael Atiyah and Isadore Singer at Harvard and M.I.T. in the fall of 1962. The index theorem was first announced in a joint paper by Atiyah and Singer \cite{AS63} in the \textit{Bulletin of the American Mathematical Society} in early 1963, using cobordism theory to obtain an expression for the index in terms of characteristic classes. In 1968, Atiyah and Singer published a second proof of the index theorem in a series of three papers in the \textit{Annals of Mathematics} (the second paper in the series is actually by Atiyah and Graeme Segal). The second proof reformulates the index theorem in terms of $K$-theory, avoiding cohomology and cobordism entirely, and providing generalizations to the equivariant $G$-index theorem and families of elliptic operators. Since 1968, the index theorem has been proved in a variety of different ways (such as the heat equation proof and the supersymmetric proof), and it has become significant to mathematical physics.

In this paper, we begin in \S\ref{sec:1} by discussing $K$-theory, the basic language of the index theorem. $K$-theory was developed primarily by Atiyah and Hirzebruch in the early 1960's, and for a general reference, the reader is referred to Atiyah's classic text \cite{Ati67}. \S\ref{sec:1} concludes with a brief review of characteristic classes and the relationship between $K$-theory and cohomology. In \S\S\ref{sec:2}--\ref{sec:3}, we will present the proof of the index theorem along the lines of \cite{AS68I} and \cite{AS68III}, digressing in \S\ref{sec:2} to present the cohomological form of the index theorem and a few simple applications. The equivariant generalization, taking into account the action of a compact Lie group $G$, is discussed in \S\ref{sec:4}. Equivariant $K$-theory was developed by Atiyah and Segal at Oxford, and a complete treatment is given by Segal in \cite{Seg68}. In \S\ref{sec:4}, we will present the $G$-index theorem and then proceed to consider the special case of homogeneous differential operators, following the paper \cite{Bott65} by Raoul Bott. In the first appendix, we will give a proof of the Bott Periodicity theorem using the index of a family of differential operators. This discussion essentially follows Atiyah's paper \cite{Ati69}, and includes a generalization of this technique to prove both the standard and equivariant versions of the Thom Isomorphism theorem. In the second appendix, we will derive the Atiyah-Hirzebruch spectral sequence without assuming any background in homological algebra.

The material presented in this text is not original, and in fact, all of the results were known before the author was born. Rather, this text is intended as a synthesis of various developments in $K$-theory and index theory during the critical period between 1963 and 1968, focusing on the works of Atiyah and his various co-authors. An attempt is made to unify the notation across the various papers in the field, and several examples are worked out in detail. In particular, the treatment of the index on odd-dimensional manifolds in \S\ref{sec:2} and the examples of clutching functions in the appendix are carried out in significantly more detail than in the available literature.

This text is intended as an introduction to the index theorem for advanced undergraduates or graudate students with little or no background in $K$-theory or index theory. However, it is assumed that the reader is familiar with the basic theory of vector bundles. If not, the reader should browse through the first few sections of \cite{MS74} or \S 1 of \cite{Ati67} before attempting to understand this text. A strong background in cohomology theory and some knowledge of characteristic classes is recommended, but it is not strictly necessary for the proof of the index theorem itself. Some familiarity with representation theory would be helpful for \S\ref{sec:4}, and the reader is referred to \cite{Ada69}.

\subsection{Definitions and Examples}

Let $X$ be a compact, smooth manifold. Given smooth complex vector bundles $E$ and $F$ over $X$, we write $\Gamma E$ and $\Gamma F$ for the spaces of $C^{\infty}$-sections (i.e., $\Gamma V$ consists of smooth maps $s : X \to V$ such that $\pi \circ s = \mathrm{Id}$, where $\pi: V\to X$ is the projection map of a vector bundle $V$) of $E$ and $F$ respectively.

\begin{definition}
	A map $P: \Gamma E \to \Gamma F$ is called a \emph{linear partial differential operator} if given local coordinates $x_{1},\ldots, x_{n}$ on $X$ it takes the form
	$$ P = \sum_{r \leq m}\alpha_{i_{1}\cdots i_{r}} \,\frac{\partial^{i_{1}+\cdots i_{r}}}{\partial x_{i_{1}}\cdots \partial x_{i_{r}}},$$
 where $\alpha_{i_{1}\cdots i_{r}}(x):E_{x}\to F_{x}$ is a linear transformation between the fibers of $E$ and $F$ depending smoothly on $x$.
\end{definition}

In other words, $P$ is locally a polynomial in the operators $\partial/\partial x_{j}$ with matrix-valued coefficients varying smoothly with $x$. Formally substituting the indeterminates $i\,\xi_{j}$ for $\partial/\partial x_{j}$, we denote the corresponding polynomial by $p(x,\xi)$, where $\xi = ( \xi_{1},\ldots,\xi_{n})$. To make clear the relationship between the polynomial and the operator, the operator can be written in the form $P = p(x,D)$, where $D = (-i\,\partial/\partial x_{1}, \ldots, -i\,\partial/\partial x_{n}$).\footnote{The reader may recognize this map  $\partial/\partial x_{j} \to i\,\xi_{j}$ as the Fourier transform. This connection with Fourier analysis will become significant later in this text.} If we remove all by the highest order terms in the polynomial $p(x,\xi)$, we obtain a homogeneous polynomial which we will denote by $\sigma(x,\xi)$ and call the \emph{symbol} of $P$.

\begin{definition}
	A linear partial differential operator is called \emph{elliptic} if its symbol $\sigma(x,\xi)$ is invertible for all $x\in X$ and $\xi = (\xi_{1}, \ldots , \xi_{n}) \neq 0 \in \R^{n}.$
\end{definition}

\begin{example}
Let $X = \R^{n}$, and let $E = F = \R^{n} \times \C$ be the trivial complex line bundle. Then an element of $\Gamma E$ or $\Gamma F$ is simply a smooth complex-valued function on $\R^{n}$. The Laplacian on $\R^{n}$ is the second order differential operator
$$\nabla^{2} = \left(\frac{\partial}{\partial x_{1}}\right)^{2} + \cdots+
\left(\frac{\partial}{\partial x_{n}}\right)^{2}.$$
Its symbol is $\sigma(x,\xi) = -\bigl( (\xi_{1})^{2}+ \cdots + (\xi_{n})^{2} \bigr) < 0$ for $\xi\neq 0$ and so the Laplacian is elliptic. More generally, if we take any operator of the form $D = f\nabla^{2}$, where $f$ is a smooth, nonzero, complex-valued function on $\R^{n}$, then $D$ is an elliptic operator as well.
\end{example}

\begin{example}
	Letting $X = \C = \R^{2}$, we set $z = x + iy$ with real coordinates $x$ and $y$. Then the Cauchy-Riemann operator from complex analysis is
	$$\frac{\partial}{\partial \bar{z}} = \frac{1}{2} \left( \frac{\partial}{\partial x} + i \,\frac{\partial}{\partial y} \right) .$$
	Recall that the Cauchy-Riemann differential equations are given by $\partial / \partial \bar{z}\,f(x,y) = 0$. The symbol of $\partial/\partial \bar{z}$ is $\sigma(x,\xi) = \frac{1}{2}\,i\, (\xi_{1} + i\xi_{2}) \neq 0$ for $\xi\neq 0$ since $\xi$ is real, and so $\partial/\partial \bar{z}$ is an elliptic operator. We will return to this operator at the end of \S\ref{sec:2}.
\end{example}

Given an elliptic operator $D : \Gamma E \to \Gamma F$ on a compact manifold $X$, the theory of partial differential operators tells us that the kernel of $D$ (the space of all solutions to $Df = 0$) and the cokernel of $D$ (the space $\Gamma F$ modulo the image of $D$) are both finite-dimensional. We can then define

\begin{definition}
The \emph{index} of an elliptic operator $D$ on a compact manifold is given by
$$\Index D = \dim\Ker D - \dim \Coker D.$$
\end{definition}
\begin{example}
Let $X$ be the circle $S^{1} = \R/2\pi\Z$, and consider the operator $d/d\theta$. Its symbol is the function $\sigma(x,\xi) = i\xi$, and so this operator is elliptic. The kernel is the space of all solutions to the differential equation $d/d\theta\,f(\theta) = 0$, which is clearly just the constant functions. To compute the cokernel, we note that any smooth function $f(\theta)$ on $S^{1}$ can be written as a Fourier series $f(\theta) = \sum_{-\infty}^{\infty}a_{n}\,e^{in\theta}$. Then, we obtain
$$\frac{d}{d\theta}\,f(\theta) = \sum_{-\infty}^{\infty}in\,a_{n}\,e^{in\theta}.$$
Noting that the constant term ($n=0$) vanishes in this series, we see that the cokernel is again just the constant functions. Hence, we have $\Index d/d\theta = 1 - 1 = 0$.

More generally, consider the elliptic operator given by $P = p(d/d\theta)$, where $p(\xi)$ is a polynomial with constant coefficients. In terms of the Fourier series expansion, it is immediate that
$$Pf(\theta) = \sum_{-\infty}^{\infty}p(in)\,a_{n}\,e^{in\theta}.$$
Letting $n_{1},\ldots,n_{r}$ be the distinct integer solutions to the equation $p(in)=0$, the kernel of $P$ consists of all functions of the form $f(\theta) = b_{1}\,e^{i n_{1}\theta} + \cdots + b_{r}\,e^{i n_{r}\theta}$, and so $\dim\Ker P = r$. Then, noting that the image of $P$ consists of all functions $f(\theta) = \sum_{-\infty}^{\infty}b_{n}\,e^{in\theta}$ with $b_{x_{1}} = \cdots = b_{x_{r}} = 0$, we see that we can identify the cokernel of $P$ with the kernel of $P$, and so we have $\dim\Coker P = r$. It follows that $\Index P = r - r = 0$. In fact, we will see later in this paper that \emph{any} linear elliptic partial differential operator on the circle $S^{1}$ has index zero.
\end{example}

As we have defined it, the index of an elliptic operator is an integer associated with the solutions of certain partial differential equations. It is the goal of this paper to express the index in terms of more easily computed topological invariants involving vector bundles, cohomology groups, and Chern classes. In the spirit of algebraic topology, our approach to this problem is summarized quite nicely by the following diagram:
$$
\xymatrix{
K(TX) \ar[r]^{\text{t-ind}} & \Z \\
\text{elliptic operators} \ar[u]^{\sigma} \ar[ur]
}
$$
The next three sections will be devoted to explaining the various facets of this diagram. We have already seen the definition of the index map from elliptic operators to the integers. \S\ref{sec:1} will be devoted to defining and studying the properties of the function $K(\cdot)$ and its relation to cohomology. In \S\ref{sec:2}, we will introduce the $K$-theoretic symbol map $\sigma$ and the topological index, $\text{t-ind}$. The actual proof of the Atiyah-Singer Index Theorem will be given in \S\ref{sec:3}, and amounts to showing that the above diagram is commutative. In the process, it will be necessary to extend our discussion to the class of elliptic \emph{pseudo}-differential operators.

\section{$K$-Theory and Cohomology}
\label{sec:1}

\subsection{Vector Bundles and $K(X)$}

Let $X$ be a compact, Hausdorff space, and let $\Vect(X)$ be the isomorphism classes of complex vector bundles over $X$. The Whitney sum $\oplus$ makes $\Vect(X)$ into an abelian semigroup with identity. We can then construct the associated group $K(X)$ of virtual vector bundles by Grothendieck's process of formally taking differences.

\begin{definition}
$K(X)$ is the set of cosets of $\Delta\Vect(X)$ in $\Vect(X) \times \Vect(X)$, where $\Delta: \Vect(X) \to \Vect(X) \times \Vect(X)$ is the diagonal map. In other words, we define $K(X)$ to be all pairs $E-F$ of vector bundles, modulo
the equivalence relation
$$E^{1} - E^{2} = F^{1} - F^{2} \Longleftrightarrow \exists G \text{ such that } E^{1}\oplus F^{2}\oplus G \cong E^{2} \oplus F^{1}\oplus G.$$
The resulting set $K(X)$ is then an abelian group with negatives given by $-(E-F) = (F-E)$.
\end{definition}

Alternately we could have defined $K(X)$ to be the free abelian group generated by $\Vect(X)$ modulo the subgroup generated by elements of the form $E + F - (E\oplus F)$. In either case, the tensor product $\otimes$ of vector bundles extends to a product on $K(X)$, giving it the structure of a commutative ring with identity.

Given a vector bundle $E$ over $X$, we denote its image in $K(X)$ by $[E]$. Note that this map $\Vect(X) \to K(X)$ need not be injective. For instance, if we consider \emph{real} vector bundles over the sphere $S^{2}\subset \R^{3}$, we recall that the normal bundle $NS^{2}$ is isomorphic to the trivial line bundle $\mathbf{1}$. In $K_{\R}(S^{2})$ we then obtain
$$[TS^{2}] = [\mathbf{3}] - [NS^{2}] = [\mathbf{3}] - [\mathbf{1}] = [\mathbf{2}].$$
while the tangent bundle $TS^{2}$ and the trivial plane bundle $[\mathbf{2}]$ are not isomorphic. This is why it was necessary to put the $G$ term in the equivalence relation defined above. So, if $[E] = [F]$, we know that there exists a vector bundle $G$ such that $E \oplus G \cong F \oplus G$. Then, there exists a vector bundle $G^{\perp}$ such that $G \oplus G^{\perp} \cong \mathbf{n}$ for some $\mathbf{n}$. It follows that $[E] = [F]$ if and only if $E$ and $F$ are stably equivalent.

Given a smooth map $f : X \to Y$ and a vector bundle $F$ over $Y$, we can construct the \emph{induced bundle} $f^{*}F$ over $X$ satisfying $(f^{*}F)_{x} = F_{f(x)}$ for all $x\in X$. This map then extends to a ring homomorphism $f^{*}:K(Y)\to K(X)$ called the \emph{induced map} of $f$. Hence, we see that $K(\cdot)$ is a contravariant functor from the category of compact spaces to the category of commutative rings with unit. In particular, if the map $i : \mathrm{pt} \to X$ is the inclusion of a point in $X$, then we obtain a map $i^{*}: K(X) \to K(\mathrm{pt}) \cong \Z$.\footnote{A vector bundle over a point is a vector space, which is given up to isomorphism by its dimension.} Defining the \emph{reduced $K$-group} by $\tilde{K}(X) = \Ker i^{*}$, we obtain the splitting $K(X) \cong \tilde{K}(X) \oplus K(\mathrm{pt})$. Viewing $\tilde{K}(X)$ as $K(X)$ modulo the trivial bundles, $\tilde{K}(X)$ then consists of classes of \emph{stable vector bundles}.

\subsection{Complexes and Compact Supports}

Now, we consider the case where $X$ is not compact, but only \emph{locally compact} (and Hausdorff).
 Letting $X^{+}$ be the one-point compactification of $X$ (in case $X$ is already compact, we define $X^{+} = X\cup \mathrm{pt}$), we can consider \emph{$K$-theory with compact supports}.

\begin{definition}
If $X$ is locally compact, then $K(X) := \tilde{K}(X^{+}) = K(X^{+}) / K(\mathrm{pt}).$
\end{definition}

Note that if $X$ is already compact, then this definition coincides with our previous one. Restricting ourself to proper maps,\footnote{A \emph{proper map} is a map for which the inverse of a compact set is compact. This property allows a proper map to extend to a continuous map on the one-point compactification of the spaces.} this definition extends $K(\cdot)$ to a contravariant function from the category of locally compact Hausdorff spaces to the category of commutative rings (note that $K(\cdot)$ is a ring \emph{with unit} if and only if $X$ is compact).\footnote{Using compact supports, $K(\cdot)$ can also be viewed as a \emph{covariant} function with respect to inclusions of open sets. Letting $i:U\to X$ be the inclusion of an open subset $U$ of $X$, we define the \emph{natural extension homomorphism} $i_{*}: K(U) \to K(X)$ to be the map
induced by $X^{+} \to X^{+}\,/\,(X^{+}- U^{+})\cong U^{+}$.}

Alternatively, we can consider complexes $E^{\bullet}$ of vector bundles over $X$ of the form
$$
\xymatrix{
0 \ar[r] & E^{0} \ar[r]^{\alpha_{0}} & E^{1} \ar[r]^{\alpha_{1}} & \cdots \ar[r]^{\alpha_{n-1}}& E^{n} \ar[r] & 0,
}
$$
where $\alpha_{k}\alpha_{k-1}=0$. Two complexes $E^{\bullet}$ and $F^{\bullet}$ over $X$ are called \emph{homotopic} if there exists a complex $G^{\bullet}$ over $X \times [0,1]$ such that $E^{\bullet}$ and $F^{\bullet}$ are isomorphic to the restrictions $G^{\bullet}|_{X\times 0}$ and $G^{\bullet}|_{X\times 1}$ respectively. The \emph{support} of this complex is the set of all points $x\in X$ where the complex restricted to the fibers at $x$ fails to be exact. For our purposes, we are interested in complexes with compact support.

\begin{claim}
Letting $C(X)$ denote the set of homotopy classes (via compactly supported homotopies) of complexes over $X$ with compact support, and letting $C_{\varnothing}(X)$ be the subset of $C(X)$ consisting of complexes with empty support, we have
$$K(X) \cong C(X) / C_{\varnothing}(X),$$
where $K(X)$ is the $K$-group with compact supports defined above.
\end{claim}

See \cite[\S2.6]{Ati67} or \cite[Appendix]{Seg68} for a proof. If $X$ is actually compact, then a complex $E^{\bullet}$ corresponds to the element
$$ \chi(E^{\bullet}) = \sum_{k}(-1)^{k}\,[E^{k}] \in K(X).$$
This alternating sum map $\chi:C(X) \to K(X)$ is known as the \emph{Euler characteristic}.

Instead of taking complexes of arbitrary length, we can work with complexes of fixed length $n > 0$, and we still obtain $K(X) \cong C^{n}(X)/C^{n}_{\varnothing}(X)$. In particular, it is often convenient to restrict ourselves to complexes of length $1$, i.e., ``pairs of vector bundles $(E,F)$ over $X$ which are isomorphic outside a compact set''. Letting $\varphi : C^{1}(X) \to C(X)$ be the obvious inclusion map, we construct a left inverse $\psi : C(X) \to C^{1}(X)$ as follows.
Given a complex $E^{\bullet}$ of arbitrary length
$$
\xymatrix{
0  \ar[r] & E^{0} \ar[r]^{\alpha_{0}} & E^{1} \ar[r]^{\alpha_{1}} & \cdots \ar[r]^{\alpha_{n-1}} & E^{n} \ar[r] & 0.
}
$$
we choose a Hermitian inner product $\langle \cdot, \cdot \rangle_{k}$ on each of the vector bundle $E^{k}$. Then, we have the adjoint maps $\alpha_{k}^{*}: E^{k+1}\to E^{k}$ defined by the property
$$\bigl\langle \alpha_{k}(v),w\bigr\rangle_{k+1} = \bigl\langle v,\alpha^{*}_{k}(w) \bigr\rangle_{k}.$$
Now, we let $\psi(E^{\bullet})$ be the corresponding complex of length one given by
$$
\xymatrix{
0 \ar[r] & \bigoplus_{i}E^{2i} \ar[r]^{\alpha} & \bigoplus_{i}E^{2i+1} \ar[r] & 0,
}
$$
where
$$\alpha = \bigoplus_{i} \alpha_{2i} + \bigoplus_{i}\alpha^{*}_{2i+1}.$$
It is easily verified that $E^{\bullet}$ and $\psi(E^{\bullet})$ have the same support. Also, note that the complex $\psi(E^{\bullet})$ is independent of our choice of inner products since all choices of inner product are homotopic to one another.

The advantage of using complexes of arbitrary length is that it provides us with a convenient way to define multiplication. Given vector bundles $E$ over $X$ and $F$ over $Y$, we define their \emph{exterior tensor product} over $X \times Y$ to be $E \tensorhat F := \pi^{*}_{X}E \otimes \pi^{*}_{Y}F$, where
$\pi_{X}$ and $\pi_{Y}$ are the coordinate projections from $X\times Y$ onto $X$ and $Y$ respectively. Taking exterior tensor products of complexes, we can extend this to a product map
$$ \tensorhat : K(X) \otimes K(Y) \to K(X\times Y),$$
where $(E^{\bullet} \oplus F^{\bullet})^{i} = \bigoplus_{j}E^{j}\otimes F^{i-j}$ and
$$ \alpha_{i}^{E^{\bullet} \tensorhat F^{\bullet} }
= \bigoplus_{j} \Bigl( \bigl(\alpha_{j}^{E^{\bullet}} \tensorhat 1_{F^{i-j}}\bigr) + (-1)^{j} \bigl( 1_{E^{j}} \tensorhat \alpha_{i-j}^{F^{\bullet}}\bigr) \Bigr).$$

For instance, given the length one complexes
$$
\xymatrix{
	0 \ar[r] & E^{0} \ar[r]^{\alpha} & E^{1} \ar[r] & 0
} \qquad \text{ and } \qquad
\xymatrix{
	0 \ar[r] & F^{0} \ar[r]^{\beta} & F^{1} \ar[r] & 0
} 
$$
over $X$ and $Y$ respectively, their exterior tensor product is the complex
$$
0 \longrightarrow E^{0} \tensorhat F^{0} \xrightarrow{\alpha\tensorhat 1 + 1 \tensorhat\beta} E^{1}\tensorhat F^{0} \oplus E^{0} \tensorhat F^{1}
\xrightarrow{-1\tensorhat\beta + \alpha\tensorhat 1}
E^{1}\tensorhat F^{1} \longrightarrow 0
$$
over $X \times Y$. Note that it is necessary to introduce the factor $-1$
so that $\alpha_{1}\alpha_{0} = 0$. By the above construction, the corresponding complex of length one is given by
$$
\xymatrix{
0 \ar[r] & E^{0}\tensorhat F^{0} \oplus E^{1} \tensorhat F^{1} \ar[r]^{\theta} &
E^{0}\tensorhat F^{1} \oplus E^{1} \tensorhat F^{0} \ar[r] & 0\,,
}
$$
where
$$\theta = \begin{pmatrix} \alpha \tensorhat 1 & -1 \tensorhat \beta^{*} \\
1 \tensorhat \beta & \alpha^{*} \tensorhat 1
\end{pmatrix}.$$

\subsection{Homogeneous Complexes}

Our primary motivation for introducing $K$-theory with compact supports is so that we can consider complexes over $V$, where $V$ is a \emph{real} vector bundle over $X$ with projection map $\pi : V \to X$ (later, we will take $V$ to be $TX$, the cotangent bundle of $X$). In this case, we can impose a homogeneous structure on the complexes. Given complex vector bundles $E$ and $F$ over $X$, we can lift them to vector bundles $\pi^{*}E$ and $\pi^{*}F$ over $V$, and we note that we have $(\pi^{*}E)_{v} = E_{\pi(v)}$ and $(\pi^{*}F)_{v} = F_{\pi(v)}.$ We say that a homomorphism $\alpha : \pi^{*}E\to\pi^{*}F$ is (positively) \emph{homogeneous of degree $m$} if for all $v\in V$ and real $\lambda > 0$, we have
$$\alpha_{\lambda v} = \lambda^{m}\alpha_{v} : E_{\pi(v)}\to F_{\pi(v)}.$$
Given a metric on $V$, we see that $\alpha$ is completely determined by its restriction to the sphere bundle $S(V)$. Suppose that $X$ is compact and $E^{\bullet}$ is a complex over $V$ given by
$$
\xymatrix{
0 \ar[r] & \pi^{*}E^{0} \ar[r]^{\alpha} & \pi^{*}E^{1} \ar[r]^{\alpha}
& \cdots \ar[r]^{\alpha} & \pi^{*}E^{n} \ar[r] & 0\,,
}
$$
where $\alpha^{2}= 0$ and $\alpha$ is homogeneous of degree $m$. If $E^{\bullet}$ is exact on $S(V)$, then the support of $E^{\bullet}$ will
be the zero section which is the image of $X$ in $V$. Thus, $E^{\bullet}$ has compact support, and so it represents an element of $K(V)$.\footnote{Note that if $m<0$ then $\alpha$ will be discontinuous on the zero section. Fortunately, we will not need the case where $m < 0$. However, if $X$ is only locally compact, then it will be necessary to take $m=0$, in which case $\alpha$ may be discontinuous on the zero section. See \cite[pp.~492--3]{AS68I}.}
We claim that we can define $K(V)$ in terms of such homogeneous complexes.

\begin{claim}
Let $\Cm(V)$ denote the set of homotopy classes of homogeneous complexes of degree $m$ over $V$, and let $\Cm_{\varnothing}(V)$ be the subset consisting of exact complexes where the homomorphism $\alpha$ is constant along each fiber of the sphere bundle $S(V)$. Then
$$\Cm(V) /\, \Cm_{\varnothing}(V) \cong C(V) / C_{\varnothing}(V) \cong K(V).\footnote{As before, we note that we may restrict our discussion to homogeneous complexes of fixed length.}$$
\end{claim}
\begin{proof}
Considering the inclusion $\Cm(V) \subset C(V)$, we would like to construct an inverse map $C(V) \to \Cm(V)$. So, given a complex $E^{\bullet}\in C(V)$ with compact support $L$
$$
\xymatrix{
	0 \ar[r] & E^{0} \ar[r]^{\alpha_{0}} & E^{1} \ar[r]^{\alpha_{1}}
	& \cdots \ar[r]^{\alpha_{n-1}} & E^{n} \ar[r] & 0\,,
}
$$
we choose a metric on $V$, and we let $D(V)$ be a disc bundle containing $L$. Then taking the restriction $F^{i} = E^{i}|_{X}$ to the zero section, we see that $E^{i} \cong \pi^{*}F^{i}$ on $D(V)$. Defining $\beta_{i}: \pi^{*}F^{i} \to \pi^{*}F^{i+1}$ to be the homogeneous map of degree $m$ that agrees with $\alpha_{i}$ on the corresponding sphere bundle $S(V)$, we thus obtain a homogeneous complex $F^{\bullet}$ on $V$
$$
\xymatrix{
	0 \ar[r] & \pi^{*}F^{0} \ar[r]^{\beta_{0}} & \pi^{*}F^{1} \ar[r]^{\beta_{1}}
	& \cdots \ar[r]^{\beta_{n-1}} & \pi^{*}F^{n} \ar[r] & 0\,.
}
$$
It is easy to see that $E^{\bullet}$ and $F^{\bullet}$ are homotopic, and so they both correspond to the same element of $K(V)$. This map $C(V) \to \Cm(V)$ is thus an isomorphism. Furthermore, if $E^{\bullet}\in C_{\varnothing}(V)$ has empty support, then we note that $E^{\bullet}$ is homotopic to a complex where $\alpha_{i}$ is constant along each fiber of $V$, and it follows that $F^{\bullet}\in \Cm_{\varnothing}(V)$. Conversely, if $F^{\bullet}\in \Cm_{\varnothing}(V)$ is a complex where $\beta_{i}$ is constant along each fiber of $S(V)$, then taking the homotopy $\alpha_{i}(v,t) = \| v\|^{tm}\beta_{i}(v / \|v\|)$, we see that $F^{\bullet}$ is homotopic to an exact complex $E^{\bullet}\in C_{\varnothing}(V)$. Thus, we have $\Cm(V) /\, \Cm_{\varnothing}(V) \cong C(V) / C_{\varnothing}(V)$. 
\end{proof}

\subsection{The Thom Isomorphism}

Let $V$ be a complex vector bundle over a compact space $X$. The exterior algebra\footnote{The \emph{nth exterior power} $\Lambda^{n}(V)$ of $V$ is the skew-symmetrization of the $n$-fold tensor product $V^{\otimes n}$.} $\Lambda^{*}(V)$ then yields a homogeneous complex $\Lambda(V)$ of degree $1$ of vector bundles over $V$ called the \emph{exterior complex}
$$
\xymatrix{
	0 \ar[r] & \pi^{*}\Lambda^{0}(V) \ar[r]^-{\alpha} &
	\pi^{*}\Lambda^{1}(V) \ar[r]^-{\alpha} & \cdots \ar[r]^-{\alpha} &
	\pi^{*}\Lambda^{n}(V) \ar[r] & 0\,.
}
$$
where $\pi: V \to X$ and $\alpha : (v,w) \mapsto (v,v\wedge w)$ for all $v\in V$ and $w\in \pi^{*}\Lambda^{k}(V)_{v} = \Lambda^{k}(V)_{\pi(v)}.$ Since $\alpha^{2}=0$, and since the complex is exact outside the zero section, we see that it defines an element $\lambda_{V}\in K(V) = \tilde{K}(V^{+}),$ where
$V^{+} = X^{V} = D(V)/S(V) = P(V\oplus \mathbf{1})/P(V)$ is the Thom space\footnote{Note that here we define the Thom space to be the one-point compactification of the total space. If $X$ is only locally compact, then we must instead take the one-point compactification of each fiber separately.} of $V$. We can now state the fundamental theorem of $K$-theory.

\begin{theorem}[Thom Isomorphism Theorem]
	If $V$ is a complex vector bundle over a compact space $X$, and $K(V)$ is viewed as a $K(X)$-module via the lifting $\pi^{*}:K(X) \to K(V)$, then the Thom map $\varphi : K(X) \to K(V)$ given by multiplication by $\lambda_{V}$ is an isomorphism.
\end{theorem}

The proof of this theorem is sketched in the appendix. If $X$ is only locally compact, then the complex $\Lambda(V)$ does not have compact support. However,
if $E^{\bullet}$ is a complex over $V$ with compact support, then the product of $\Lambda(V)$ and $E^{\bullet}$ does have compact support, and the Thom Isomorphism Theorem still holds. Letting $s : X \to V$ be the zero section, we note that if $X$ is compact, then
$$ s^{*}\lambda_{V} = s^{*}\circ \varphi(1) =
\sum_{k=0}^{n}(-1)^{k}\Lambda^{k}(V).$$
Noting that $\Lambda^{*}(V \oplus W) \cong \Lambda^{*}(V) \otimes \Lambda^{*}(W)$, we see that $\lambda_{V\oplus W} = \lambda_{V} \cdot \lambda_{W}$, and hence the Thom isomorphism is transitive. In other words, the isomorphism $K(X) \to K(V\oplus W)$ is identical to the composition $K(X) \to K(V) \to K(V \oplus W).$

\begin{example}
If $X$ is a point and $V = \C^{n}$, then we have the Thom isomorphism
$$\varphi : K(\mathrm{pt}) \xrightarrow{\cong} K(\C^{n}) \cong \tilde{K}(S^{2n}),$$
and thus $\tilde{K}(S^{2n})$ is isomorphic to $\Z$ and is generated by the \emph{Bott class} $\lambda_{\C^{n}} = \lambda_{n}.$ It then follows that $K(S^{2n}) \cong \Z^{2}$.
\end{example}

\subsection{A Periodic Cohomology Theory}

One of the most interesting and useful features of $K$-theory is that it admits the properties characteristic of a generalized cohomology theory. To see this, we will first require a few definitions.

\begin{definition}
The \emph{reduced suspension} of a space $X$ with base point is given by
$$SX = S^{1}\wedge X = (S^{1} \times X) / ( S^{1} \vee X),$$
and the $n$-th iterated suspension of $X$ is given by
$$S^{n}X = \underbrace{SS\cdots S}_{n} X = S^{n} \wedge X,$$
where $X \vee Y = (X \times \mathrm{pt}) \cup ( \mathrm{pt} \times Y)$
is the \emph{one point union} of two spaces with base point, and $X \wedge Y = X \times Y \,/\, X \vee Y$ is the \emph{smash product}.
\end{definition}

\begin{definition}
For $n\geq 0$, we define degree shifts by\footnote{Using $K$-theory with compact supports, an equivalent definition is $K^{-n}(X) := K(\R^{n}\times X)$.}
\begin{align*}
	K^{-n}(X,Y) &:= \tilde{K}(S^{n}(X/Y)), \\
	K^{-n}(X) &:= K^{-n}(X,\varnothing) = \tilde{K}(S^{n}X^{+}).\end{align*}
\end{definition}

\begin{corollary}[Bott Periodicity]
The map given by
$$K^{-n}(X) \xrightarrow{\lambda_{1}\tensorhat\cdot} K^{-n-2}(X),$$
where $\lambda_{1}$ is the Bott class generating $\tilde{K}(S^{2}) = K^{-2}(\mathrm{pt}),$ is an isomorphism for all $n \geq 0$.
\end{corollary}

\begin{proof}
Letting $V = \C^{m}\times X$ be the trivial bundle, we have the Thom isomorphism
\begin{equation*}\begin{split}
K(X) \xrightarrow{\cong} K(\C^{m}\times X) &= \tilde{K}\bigl( (\C^{m}\times X)^{+} \bigr) \\ &= \tilde{K}(S^{2m} \wedge  X^{+}) = K^{-2m}(X)
\end{split}\end{equation*}
given by
$$x\mapsto \lambda_{X\times \C^{m}} \cdot \pi^{*}x = \pi^{*}_{\C^{m}}\lambda_{m} \cdot \pi^{*}_{X}x = \lambda_{m} \tensorhat x,$$
where we view the exterior tensor product as a map $\tensorhat : \tilde{K}(X) \otimes \tilde{K}(Y) \to \tilde{K}(X\wedge Y)$. Noting that $\lambda_{m} = (\lambda_{1})^{m}$, the desired result follows immediately.
\end{proof}

By Bott periodicity, we see that all of the even dimensional groups are isomorphic to $K^{0}(X) = K(X)$ and all of the odd dimensional groups are isomorphic to $K^{1}(X) = \tilde{K}(SX^{+}).$ We thus obtain a $\Z_{2}$-graded ring $K^{*}(X) = K^{0}(X) \oplus K^{1}(X)$, which gives us a periodic cohomology theory with the following exact ``ring'' sequence:
$$
\xymatrix{
	K^{0}(X,Y) \ar[r] & K^{0}(X) \ar[r] & K^{0}(Y) \ar[d] \\
	K^{1}(Y) \ar[u] & K^{1}(X) \ar[l] & K^{1}(X,Y) \ar[l]
}
$$

Using this exact sequence, we can obtain results analogous to those of integral cohomology. For example, we leave the proof of the following proposition to the reader (or see \cite[\S2.5.2]{Ati67}).

\begin{proposition}
	If $X$ is a finite $CW$-complex with cells in only even dimensions, then $K^{0}(X) \cong \Z^{r},$ where $r$ is the number of cells, and $K^{1}(X) = 0$. Hence $K^{*}(X) \cong \Z^{r}$.
\end{proposition}

So, recalling that complex projective space $\C P^{n}$ has exactly one cell in each of the dimensions $0,2,\ldots,2n$, we see that $K^{*}(\C P^{n})\cong \Z^{n+1}$. We also obtain $K(S^{2n}) \cong \Z^{2}$ as before. We now present the analog of the K\"unneth formula in $K$-theory (see \cite[\S2.7.15]{Ati67}).

\begin{theorem}[K\"unneth Theorem]
Let $X$ and $Y$ be finite CW-complexes. Then we have the natural exact sequence (with indices in $\Z_{2}$)
$$
0 \to \sum_{i+j=k}K^{i}(X) \otimes K^{j}(Y) \to K^{k}(X\times Y) \to
\sum_{i+j=k+1}\mathrm{Tor}\bigl(K^{i}(X), K^{j}(Y) \bigr) \to 0\,.
$$
\end{theorem}

\begin{example}
Since all complex vector bundles over the circle $S^{1}$ are trivial, we have $K(S^{1})\cong \Z$. So, for the case of a point, we obtain $K^{0}(\mathrm{pt}) \cong \Z$ and $K^{1}(\mathrm{pt}) = \tilde{K}(S^{1}) = 0$. On the other hand,
for the circle $S^{1}$, we see that $K^{0}(S^{1}) \cong \Z$ and $K^{1}(S^{1}) = \tilde{K}(S^{2}) \cong \Z$. Then, by the K\"unneth Theorem, we obtain $K(X \times S^{1}) \cong K^{0}(X) \oplus K^{1}(X) = K^{*}(X)$.\footnote{We could also take $K(X\times S^{1} ) \cong K(X \wedge S^{1}) \oplus K(X) \oplus K(S^{1}) \oplus K(\mathrm{pt}) \cong \tilde{K}(SX^{+}) \oplus K(X).$}
\end{example}

\subsection{The Chern Character}

In order to discuss the relationship between $K$-theory and cohomology, it is necessary to introduce the concept of characteristic classes. Here, we will be using \v{C}ech cohomology, taking cohomology with compact supports when dealing with locally compact spaces. Given a complex $n$-plane bundle $E$ over a base space $X$, we associate with it the \emph{Chern classes} $c_{i}(E) \in H^{2i}(X;\Z)$ for $i= 0,\ldots,n$, where $c_{0}(E) = 1$, and we define the \emph{total Chern class} to be the formal sum
$$c(E) = 1 + c_{1}(E) + \cdots + c_{n}(E).$$
For a construction of the Chern classes, see \cite[\S 14]{MS74} or \cite[\S 20]{BT82}. For our purposes, we need only know that the Chern classes satisfy the following four important properties:\footnote{In order to eliminate the degenerate case where all of the Chern classes are zero, we me also require that $c_{1}(\text{Hopf bundle})$ be the canonical generator of $H^{2}(\C P^{n})$ for $n\geq 1$. This fact together with the naturality property  and the product formula are sufficient to completely characterize the Chern classes.}

\begin{property}[Naturality]
If $f: X \to Y$ is covered by a bundle map $E\to F$ between complex $n$-plane bundles $E$ over $X$ and $F$ over $Y$, then $c(E) = f^{*}c(F).$
\end{property}
\begin{property}[Product Formula]
If $E$ and $F$ are complex $m$-plane and $n$-plane bundles respectively over the same base space $X$, then the total Chern class satisfies the formula
$$ c(E\oplus F) = c(E) \smallsmile c(F) \text{, i.e., }
   c_{i}(E\oplus F) = \sum_{k=0}^{i}c_{k}(E) \smallsmile c_{i-k}(F) \text{ for all }i \geq 0.$$
\end{property}
\begin{property}
If we denote by $\overline{E}$ the conjugate bundle of a complex vector bundle $E$ (so $\overline{E}$ has the same underlying real space as $E$, but has the opposite complex structure), then $c_{k}(\overline{E}) = (-1)^{k}c_{k}(E)$, or in terms of the total Chern class,
$$ c(\overline{E}) = 1 - c_{1}(E) + c_{2}(E) - \cdots \pm c_{n}(E).$$
\end{property}
\begin{property}
If $\mathbf{n}$ is the trivial complex $n$-plane bundle, then $c(\mathbf{n}) = 1$.
\end{property}

Suppose $E$ is a complex vector bundle over $X$ that splits as the direct sum of complex line bundles $E = L_{1}\oplus \cdots L_{n}$. Then, letting $x_{i} = c_{1}(L_{i}) \in H^{2}(X;\Z)$, the product formula gives us
$$ c(E) = \prod_{i}c(L_{i}) = \prod_{i}(1+x_{i}) = 1 + \sum_{i}x_{i} + \sum_{i< j}x_{i}x_{j} + \cdots,$$
and so we see that the Chern class $c_{k}(E)$ is the $k$-th elementary symmetric function of the $x_{i}$. It follows that any symmetric function of the $x_{i}$ can be expressed as a polynomial in the Chern classes. Although not every vector bundle decomposes in this fashion,\footnote{The simplest example of a vector bundle that does not split is $T\C P^{2}$.} we do have the following important theorem.

\begin{theorem}[Splitting Principle]
If $E$ is a complex vector bundle over $X$, then there exists a space $F(E)$ and a map $\pi:F(E) \to X$ such that $\pi^{*}:H^{*}(X) \to H^{*}(F(E))$ is injective and $\pi^{*}E$ splits as the direct sum $E = L_{1}\oplus \cdots \oplus L_{n}$ of complex line bundles.
\end{theorem}

More precisely, we take $F(E)$ to be the flag bundle of $E$. For a proof, see \cite[\S 2.7.1]{Ati67} or \cite[\S 21]{BT82}. The splitting principle, when combined with the naturality of the Chern classes, tells us that in computations involving only the Chern classes, any vector bundle can be treated as if it were the direct sum of line bundles. In particular, we can define:

\begin{definition}
The \emph{Chern character} of a complex $n$-plane bundle $E$ is given by
\begin{equation*}\begin{split}
	\ch(E) = \sum_{i}e^{x_{i}} &= n + \sum_{i}x_{i} + \frac{1}{2!}\sum_{i}x_{i}^{2} + \cdots \\
	&= n + c_{1} + \frac{1}{2} ( c_{1}^{2} - 2c_{2} ) + \cdots .
\end{split}\end{equation*}
\end{definition}

For line bundles $L_{1}$ and $L_{2}$, we note that $c_{1}(L_{1}\otimes L_{2}) = c_{1}(L_{1}) + c_{1}(L_{2}),$ so we have
$$\ch(L_{1}\otimes L_{2}) = e^{x_{1}+x_{2}} = e^{x_{1}}e^{x_{2}} = \ch(L_{1}) \smallsmile \ch(L_{2}).$$
It follows that for any complex vector bundles $E$ and $F$ we obtain
\begin{align*}
	\ch(E\oplus F) &= \ch(E) + \ch(F), \\
	\ch(E \otimes F) &= \ch(E) \smallsmile \ch(F).
\end{align*}
Extending the Chern character to $K(X)$, we thus obtain a natural ring homomorphism
$$\ch : K(X) \to H^{\mathrm{even}}(X;\Q)$$
from $K(X)$ to the even dimensional rational cohomology of $X$. Furthermore, composing $\ch : \tilde{K}(SX) \to \tilde{H}^{\rm{even}}(SX^{+};\Q)$ with $\alpha^{-1}$, where $\alpha : H^{\rm{odd}}(X;\Q) \to \tilde{H}^{\rm{even}}(SX^{+};\Q)$ is the suspension isomorphism given by multiplying by the canonical generator of $H^{1}(S^{1};\Z)$, we can extend the Chern character even further to obtain a natural ring homomorphism
$$\ch: K^{1}(X) \to H^{\rm{odd}}(X;\Q).$$
Thus, we may view the Chern character as a natural transformation of functors $K^{*}(\cdot) \to H^{*}(\cdot;\Q)$, and if we eliminate torsion by using rational $K$-theory, we obtain

\begin{theorem}
If $X$ is a finite CW-complex, then the Chern character
$$\ch:K^{*}(X)\otimes \Q \to H^{*}(X;\Q)$$
is a natural isomorphism between rational $K$-theory and rational cohomology.
\end{theorem}

For a spectral sequence proof, see \cite[\S 2]{Ati-Hir}. The special case of even dimensional spheres is discussed in \cite[\S 3]{Hir72} in the context of the proof of the Bott periodicity theorem.

\subsection{The Cohomology Thom Isomorphism}

Next, suppose that $E$ is an oriented \emph{real} $n$-plane bundle over $X$. We know that there exists a unique generator $u\in H^{n}(E,E_{0};\Z) = H^{n}_{c}(E;\Z)$, where $H^{*}_{c}(\cdot)$ is cohomology with compact supports and $E_{0}$ is the deleted space obtained by removing the zero section of $E$, such that the restriction of $u$ to each fiber of $E$ induces the preferred orientation. This generator $u$ is known as the \emph{cohomology Thom class} of $E$. Then, letting $i : X \to E$ be the zero section, we define the \emph{Euler class} of $E$ by its restriction to $X$:
$$e(X) = i^{*}u \in H^{n}(X;\Z).$$
If $E$ is a complex $n$-plane bundle over $X$, then the complex structure of $E$ induces an orientation on the underlying real $2n$-plane bundle $E_{\R}$. The top Chern class is then defined by $c_{n}(E) = e(E_{\R})$. We now state the cohomology Thom isomorphism theorem.

\begin{theorem}[Thom Isomorphism]
The Thom homomorphism
$$\psi : H^{k}(X) \to H^{n+k}_{c}(E),$$
given by $\psi : x \mapsto \pi^{*}x \smallsmile u$, where $\pi : E \to X$ is the projection, is an isomorphism.\footnote{Unlike in $K$-theory, the cohomology Thom isomorphism theorem applies when $E$ is a \emph{real} bundle.}
\end{theorem}

For a proof, see \cite[\S 10]{MS74}. Using the Chern character, we can compare the Thom isomorphisms of $K$-theory and cohomology. Letting $E$ be a complex $n$-plane bundle over a compact space $X$, we obtain the following diagram:
$$
\xymatrix{
	K(X) \ar[r]^{\varphi} \ar[d]^{\ch} & K(E) \ar[d]^{\ch} \\
	H^{*}(X;\Q) \ar[r]^{\psi} & H_{c}^{*}(E;\Q)
}
$$
where the vertical maps are the corresponding Chern characters and the horizontal maps are the Thom isomorphisms $\varphi : x \mapsto \pi^{*}x \cdot \lambda_{E}$ and $\psi : x \mapsto \pi^{*}x \smallsmile u$. Since this diagram does not commute, it is necessary to introduce a correction factor. Precisely, we set
$$\psi\bigl( \mu(E) \smallsmile \ch(x) \bigr) = \ch \bigl( \varphi(x) \bigr),$$
where
$$\mu(E) = \psi^{-1} \circ \ch \circ \varphi(1) = \psi^{-1}\ch(\lambda_{E}) \in H^{*}(X;\Q).$$
It remains to calculate this cohomology class $\mu(E)$ explicitly.

Letting $i : X \to E$ be the zero section, we recall that $i^{*}\lambda_{E} = \sum_{i}(-1)^{i}\Lambda^{i}(E)$, and we note that for all cohomology classes $x \in H^{*}(X;\Q)$ we have $i^{*}\psi(x) = x \smallsmile e(E_{\R}).$ In particular, taking $x =\mu(E)$, it follows by the naturality of the Chern classes that
\begin{equation*}\begin{split}
\mu(E) \smallsmile e(E_{\R}) &= i^{*}\psi\bigl( \mu(E) \bigr) = i^{*}\ch(\lambda_{E}) \\
&= \ch( i^{*}\lambda_{E} ) = \ch \left( \sum_{i=0}^{n}(-1)^{i}\Lambda^{i}(E) \right).
\end{split}\end{equation*}
Then, invoking the splitting principle, we can assume that $E$ decomposes as the direct sum $E = L_{1}\oplus \cdots \oplus L_{n}$. Noting that $\Lambda^{i}(E) \cong \bigoplus_{j_{1}<\cdots < j_{i}}L_{j_{1}}\otimes \cdots \otimes L_{j_{i}}$ we obtain
\begin{equation*}\begin{split}
	\mu(E) \smallsmile e(E_{\R}) &= \ch \left( \sum_{i=0}^{n}(-1)^{i}\Lambda^{i}(E) \right) \\
	&= \ch\left( \prod_{i=1}^{n}(1-L_{i}) \right) = \prod_{i=1}^{n}(1-e^{x_{i}}),
\end{split}\end{equation*}
where $x_{i} = c_{1}(L_{i}) \in H^{2}(X;\Z)$ for $i=1,\ldots,n$. Since $e(E_{\R}) = c_{n}(E) = \prod_{i}x_{i}$ we conclude
$$\mu(E) = \prod_{i=1}^{n}\frac{1-e^{x_{i}}}{x_{i}}.$$

To simplify this result even further, we introduce yet another characteristic class.

\begin{definition}
The \emph{Todd class} of a complex vector bundle $E$ is given by
$$\td(E) = \prod_{i=1}^{n}\frac{x_{i}}{1-e^{-x_{i}}} = 1 + \frac{1}{2}\, c_{1}+ \frac{1}{12}\,(c_{2}+ c_{1}^{2}) + \cdots .$$
\end{definition}

Note that by its definition, the Todd class satisfies $\td(E\oplus F) = \td(E) \smallsmile \td(F)$ for any two vector bundles $E$ and $F$, and also $\td(\mathbf{n}) = 1$ for the trivial complex $n$-plane bundle $\mathbf{n}$.
We then see that, in terms of the Todd class, the correction factor is given by
$$\mu(E) = (-1)^{n}\td(\overline{E})^{-1},$$
noting that the Todd class is invertible since it has constant term $1$.

\section{Statement of the Index Theorem}
\label{sec:2}

\subsection{The Symbol Map and Elliptic Complexes}

Suppose we have an elliptic operator $D : \Gamma E \to \Gamma F$ on a compact space $X$. Then for each $x \in X$ and $\xi\in TX_{x}$, the symbol $\sigma(x,\xi)$ of $D$ is a homomorphism $E_{x}\to F_{x}$ varying smoothly with $x$ and $\xi$. It thus defines a homogeneous complex
$$
\xymatrix{
	0 \ar[r] & \pi^{*}E \ar[r]^{\alpha} & \pi^{*}F \ar[r] & 0
}
$$
over $TX$.\footnote{In this section and the next, we will use $TX$ to denote the \emph{cotangent} bundle of $X$. Note that if we choose a metric on $X$ then we can identify the tangent and cotangent bundles of $X$.}
Since $D$ is elliptic, we know that $\sigma(x,\xi)$ is invertible for $\xi\neq 0$, and so this complex is exact outside the zero section. It follows that this complex defines an element of $K(TX)$, called the \emph{symbol} of $D$ and denoted by $\sigma(D)$.

This construction can be extended with little difficulty to complexes of differential operators. Given a sequence of vector bundles $E^{0},\ldots,E^{n}$ over a compact space $X$ and partial differential operators $d_{i} : \Gamma E^{i}\to \Gamma E^{i+1},$ we say that the complex
$$
\xymatrix{
	0 \ar[r] & \Gamma E^{0} \ar[r]^-{d_{0}} & \Gamma E^{1} \ar[r]^-{d_{1}}
	& \cdots \ar[r]^-{d_{n-1}} & \Gamma E^{n} \ar[r] & 0
}
$$
is an \emph{elliptic complex} if $d_{i+1}d_{i} = 0$ and the corresponding symbol complex
$$
\xymatrix{
	0 \ar[r] & \pi^{*} E^{0} \ar[r]^-{\sigma_{0}} & \pi^{*} E^{1} \ar[r]^-{\sigma_{1}}
	& \cdots \ar[r]^-{\sigma_{n-1}} & \pi^{*} E^{n} \ar[r] & 0
}
$$
over $TX$ is exact outside the zero section. As before, this complex gives us an element $\sigma(D^{\bullet})\in K(TX)$ called the \emph{symbol} of the elliptic complex $D^{\bullet}$. In the case of elliptic complexes, the analog of the index is the \emph{Euler characteristic} defined by
$$ \chi (D^{\bullet}) = \sum_{i=0}^{n}(-1)^{i}\dim H^{i}(D^{\bullet}),$$
where $H^{i}(D^{\bullet}) = \Ker d_{i} / \im d_{i-1}$ is the cohomology of the complex. Note that in the case where the complex $D^{\bullet}$ is just a single operator, the Euler characteristic reduces to our original definition of the index. Furthermore, if we construct the corresponding complex of length one as in \S \ref{sec:1}, we obtain an elliptic operator
$$D : \Gamma \Bigl( {\bigoplus}_{i} E^{2i} \Bigr) \to \Gamma \Bigl( {\bigoplus}_{i} E^{2i+1} \Bigr)$$
given by $D = \bigoplus_{i}d_{2i} + \bigoplus_{i}d^{*}_{2i+1},$ where $d^{*}$ denotes the adjoint with respect to some metric. We then note that $\sigma(D) = \sigma(D^{\bullet})\in K(TX)$ and $\Index D = \chi(D^{\bullet})$, and so the problem is reduced to that of a single elliptic operator.

\begin{example}
Consider the (complexified) \emph{de\,Rham complex}
$$
\xymatrix{
	0 \ar[r] & \Gamma\Lambda^{0}(TX\otimes\C) \ar[r]^-{d} & \Gamma\Lambda^{1}(TX\otimes\C) \ar[r]^-{d}
	& \cdots \ar[r]^-{d} & \Gamma\Lambda^{n}(TX\otimes\C) \ar[r] & 0\,,
}
$$
where $\Gamma\Lambda^{i}(TX)$ is the space of differential forms of degree $i$ on $X$, and $d$ is the exterior derivative. In this case, the cohomology of the complex is just the de\,Rham cohomology.
$H^{i}(D^{\bullet}) = H^{i}_{\mathrm{dR}}(X) \otimes \C$, and so the Euler characteristic of the de\,Rham complex
$$\chi(D^{\bullet}) = \sum_{i}\dim H^{i}_{\mathrm{dR}}(X)
= \sum_{i}(-1)^{i}\dim H^{i}(X;\R) = \chi(X)$$
is equal to the Euler characteristic of the manifold $X$. The corresponding element $\sigma(D^{\bullet})\in K(TX)$ is called the \emph{de\,Rham symbol} and is denoted by $\rho_{X}$.
\end{example}

\subsection{Construction of the Topological Index}

Since the symbol provides us with a map from elliptic operators on $X$ into $K(TX)$, we would like to compose it with a suitable map $K(TX)\to\Z$. One such map that immediately springs to mind is the map induced by the inclusion $i : \mathrm{pt} \to TX$ of a point. However, recalling our definition $K(TX) = \tilde{K}(TX^{+})$, we observe much to our disappointment that the map $i^{*}:K(TX) \to K(\pt) \cong \Z$ is identically zero. In fact, noting that $i$ extends to a map $i : \pt \to TX^{+}$, this follows immediately since we defined the reduced $K$-group $\tilde{K}(TX^{+})$ to be $\Ker i^{*}$. So, we must look deeper for our desired map.

Instead of considering the inclusion of a base point in $X$, it is necessary to think globally. We recall that if $X$ is a compact manifold, then $X$ can be embedded in $\R^{n}$ for suitably large $n$.\footnote{This can be shown using a simple construction involving partitions of unity.} Given any such embedding $X \subset \R^{n}$, we can extend it to an embedding $TX \subset T\R^{n}$ of the tangent bundle. Since $K(T\R^{n}) \cong \Z$ by the Thom isomorphism theorem, we would be satisfied if we were to construct a map $K(TX) \to K(T\R^{n})$.

More generally, if $i : X \to Y$ is the inclusion of a compact submanifold $X$ in $Y$, we can extend it to an inclusion $i : TX \to TY$ of the tangent bundles. We then construct the map
$$i_{!} : K(TX) \to K(TY)$$
as follows (note that this map is functorial with respect to inclusions). Given a metric on $Y$, we let $N$ be a tubular neighborhood of $X$ in $Y$ diffeomorphic to the normal bundle of $X$ in $Y$. Extending to tangent bundles, we see that $TN$ is also a tubular neighborhood of $TX$ in $TY$ diffeomorphic to the normal bundle of $TX$ in $TY$. For a proof of the tubular neighborhood theorem, see \cite[\S 11]{MS74}. Since $TN$ is a \emph{real} vector bundle over $TX$, we would like to give it a complex structure so that we can apply the Thom isomorphism.

Considering an arbitrary vector bundle $E$ over $X$, we recall that $E$ is locally homeomorphic to $U \times \R^{m}$, where $U\subset X$ is homeomorphic to an open region of $\R^{n}$. Then, the tangent space to $E$ at $(x,\xi)\in E$ is given by $TE_{(x,\xi)} = TU_{x}\times T\R^{m}_{\xi}\cong TX_{x}\oplus E_{x}$, and it follows that the tangent bundle $TE$ over $E$ admits the decomposition
$$TE = \pi^{*}TX \oplus \pi^{*}E,$$
where $\pi : E\to X$ is the projection. In particular, taking $E = TX$ and $E = TY|_{X}$, we obtain
\begin{align*}
	T(TX) &\cong \pi^{*}TX \oplus \pi^{*}TX, \\
	T(TY)|_{TX} &\cong \pi^{*}(TY|_{X}) \oplus \pi^{*}(TY|X),
\end{align*}
where $T(TX)$ and $T(TY)|_{TX}$ are both bundles over $TX$ and $\pi : TX \to X$ is the projection. Then, since $N$ and $TN$ are the normal bundles of $X$ and $TX$ in $Y$ and $TY$ respectively, we have $N\oplus TX \cong TY|_{X}$ and $TN\oplus T(TX) \cong T(TY)|_{TX}$. It follows that $TN$ decomposes as
$$TN \cong \pi^{*}N \oplus \pi^{*}N.$$
We are thus able to impose on $TN$ the structure of a complex vector bundle:
$$TN \cong \pi^{*}N \oplus i\,\pi^{*}N \cong (\pi^{*}N) \otimes_{\R} \C \cong \pi^{*}( N \otimes_{\R} \C) .$$

So, composing the Thom isomorphism $\varphi : K(TX) \to K(TN)$ with the natural extension homomorphism $h : K(TN) \to K(TY),$ we obtain our desired map
$$i_{!} : K(TX) \xrightarrow{\varphi} K(TN) \xrightarrow{h} K(TY).$$
We note that this map is independent of the choice of neighborhood $N$, and that it is functorial (i.e., $(j\circ i)_{!} = j_{!}\circ i_{!}$) since the Thom isomorphism is transitive. Also, we have
$$i^{*} \circ i_{!} : x \mapsto i^{*}\lambda_{TN}\cdot x = 
\left( \sum_{i=0}^{n}(-1)^{i}\Lambda^{i}(N\otimes_{\R}\C) \right) \cdot x.$$

\begin{definition}
If $X$ is a compact manifold, choose an embedding $i:X \to \R^{m}$ for large enough $m$, and let $j:P \to \R^{m}$ be the inclusion of the origin. We note that $j_{!}:K(TP)\to K(T\R^{m})$ is simply the Thom isomorphism $\varphi : K(\pt) \to K(\C^{m})$, and so we can define the \emph{topological index} to be the composition
$$\tind : K(TX) \xrightarrow{i_{!}} K(T\R^{m}) \xrightarrow{j_{!}^{-1}} K(TP)\cong \Z.$$
\end{definition}

Suppose that we have two embeddings $i: X \to \R^{m}$ and $i' : X \to \R^{n}$. Considering the diagonal embedding $i\oplus i' : X \to \R^{m+n}$, we see that $i \oplus i'$ is isotopic to the embeddings $i\oplus 0$ and $0\oplus i'$. Letting $\lambda_{n}\in K(\C^{n})$ denote the Thom class of $\C^{n} \cong T\R^{n}$, we note that $(i\oplus 0)_{!}(x) = i_{!}(x) \cdot \lambda_{n}$ and $(0\oplus i')_{!}(x) = \lambda_{m}\cdot i'_{!}(x).$ We then see by the transitivity of the Thom isomorphism that $i\oplus 0$ and $i$ determine the same $\tind$, as do $0 \oplus i'$ and $i'$. It follows that the topological index does not depend on our choice of embedding.

Now that we have constructed our desired map $\tind:K(TX)\to \Z$, we can compose it with the symbol map $\sigma : \text{(elliptic operators)} \to K(TX)$. We are finally prepared to state the main result of this paper.

\begin{theorem}[Atiyah-Singer Index Theorem]
If $D$ is an elliptic operator on a compact manifold then
$$\Index D = \tind \sigma(D).$$
\end{theorem}

\subsection{The Cohomological Form of the Index Theorem}

Before we proceed with the $K$-theory proof of the Atiyah-Singer Index Theorem, we will briefly digress to reformulate the theorem in terms of characteristic classes and cohomology. It was in this form that the theorem was originally proved by Atiyah and Singer in 1963 using techniques from cobordism theory. Although the $K$-theoretic formulation may offer a more elegant and general proof, the cohomological form lends itself towards direct computation of examples. To write the topological index in terms of cohomology, we let $X$ be a compact $n$-dimensional manifold embedded in $\R^{k}$ and consider the following non-commutative diagram:
$$
\xymatrix{
K(TX) \ar[r]^{\varphi} \ar[d] & K(TN) \ar[r]^{h} \ar[d] & K(T\R^{k}) \ar[d] & K(TP)\cong \Z \ar[l]_{\varphi'} \ar[d] \\
H^{*}_{c}(TX) \ar[r]^{\psi}& H^{*}_{c}(TN) \ar[r]^{k} & H^{*}_{c}(T\R^{k}) & H^{*}_{c}(TP) \cong \Q \ar[l]_{\psi'}
}
$$
where $\varphi$, $\psi$, $\varphi'$, and $\psi'$ are the various Thom isomorphisms, $h$ and $k$ are the extension homomorphisms, and the vertical map is the Chern character $\ch : K(\cdot) \to H^{\mathrm{even}}_{c}(\cdot;\Q).$
Although this diagram does not commute, we recall from \S\ref{sec:1} that for a complex vector bundle $E$ over $X$, the Thom isomorphisms of $K$-theory and cohomology are related by
$$\psi^{-1}\circ \ch \circ \varphi : x \mapsto (-1)^{n}\td(\overline{E})^{-1}\smallsmile \ch(x).$$
In our case, we have $\td(\overline{T\R}^{k}) = 1$ since $\overline{T\R}^{k}$ is clearly a trivial bundle over the point $TP$. since $TN \cong \pi^{*}N \otimes_{\R}\C$, we see that $TN \cong \overline{TN}$. Noting that $T(TX)\oplus TN = T(T\R^{k})|_{TX}$ is a trivial bundle over $TX$ and that
$T(TX)\cong \pi^{*}TX \oplus \pi^{*}TX \cong \pi^{*}TX \otimes_{\R}\C$, we obtain
$$\td(\overline{TN})^{-1} = \td(TN)^{-1} = \td(T(TX)) = \pi^{*}\td(TX\otimes_{\R}\C) = \pi^{*}\,\Im(X),$$
where the last equality provides the definition of the \emph{index class} $\Im(X)$ of a manifold $X$. Hence, the class $\td(\overline{TN})^{-1} = \pi^{*}\,\Im(X)$ is independent of the embedding. We thus have 
\begin{align*}
	\psi'^{-1}\circ \ch \circ \varphi' &: x\mapsto (-1)^{k}\td(\overline{T\R}^{k})^{-1} \smallsmile \ch(x) = (-1)^{k}\ch(x), \\
	\psi^{-1}\circ \ch \circ \varphi &: x\mapsto (-1)^{n-k}\td(\overline{TN})^{-1}\smallsmile \ch(x) = (-1)^{n-k}\,\Im(X) \smallsmile \ch(x).
\end{align*}
Recalling that $\tind = \varphi'^{-1} \circ h \circ \varphi$, by an extended diagram chase we obtain
$$\ch \circ h \circ \varphi(x) = \ch\circ \varphi'(\tind x)
= (-1)^{k}\psi'(\tind x) = (-1)^{k}(\tind x)\psi'(1),$$
where $\psi'(1)$ is the canonical generator of $H^{*}_{c}(T\R^{k}).$ Then letting $[T\R^{k}]$ be the fundamental homology class of $T\R^{k}$, we have $\psi'(1)[T\R^{n}] = 1$, and thus we calculate
\begin{equation*}\begin{split}
	\tind x &= \bigl( \tind x \smallsmile \psi'(1) \bigr) [T\R^{n}] 
	= (-1)^{k}\bigl(\ch\circ h \circ \varphi(x) \bigr) [T\R^{n}] \\
	&= (-1)^{k}\bigl( \ch \circ \varphi(x) \bigr) [TN]
	= (1)^{n} \bigl( \psi(\ch(x) \smallsmile \Im(X))\bigr) [TN] \\
	&= (-1)^{n}\bigl( \ch(x) \smallsmile \Im(X) \bigr) [TX].
\end{split}\end{equation*}
We can now state a topological version of the Atiyah-Singer Index Theorem.

\begin{theorem}[Index Theorem A]
If $D$ is an elliptic operator on a compact $n$-dimensional manifold $X$ then
$$\Index D = (-1)^{n}\bigl( \ch(\sigma(D)) \smallsmile \Im(X) \bigr) [TX],$$
where $\Im(X) = \td(TX \otimes_{\R} \C) \in H^{*}(X;\Q)$ is the index class of $X$.
\end{theorem}

To obtain the precise statement from the 1963 paper on the index theorem, we note that if $X$ is an oriented manifold, then we have $\psi(u)[TX] = (-1)^{n(n-1)/2}u[X]$ for each $u\in H^{n}(X;\Z)$, where $\psi:H^{*}(X) \to H^{*}(TX)$ is the Thom isomorphism. It is necessary to introduce the sign $(-1)^{n(n-1)/2}$ because of our choice of orientation on $X$. We then obtain the original index theorem:

\begin{theorem}[Index Theorem B]
If $D$ is an elliptic operator on a compact oriented $n$-dimensional manifold $X$ then
$$\Index D = (-1)^{n(n+1)/2}\bigl( \psi^{-1}\circ \ch(\sigma(D)) \smallsmile \Im(X) \bigr) [X],$$
where $\Im(X) = \td(TX\otimes_{\R}\C)\in H^{*}(X;\Q)$ is the index class of $X$.
\end{theorem}

In order to facilitate the computation of the topological index for specific elliptic operators or complexes, we now present a special case of the Atiyah-Singer Index Theorem. Under the appropriate circumstances, we can simplify the $\psi^{-1}\circ \ch(\sigma(D))$ term in the statement of the theorem. Letting $\psi: H^{*}(X) \to H^{*}(TX)$ be the Thom isomorphism and $i^{*} : H^{*}(TX) \to H^{*}(X)$ be the map induced by the zero section $X\to TX$, we recall from \S\ref{sec:1} that $i^{*}\psi(x) = x\smallsmile e(X)$, where $e(X) = e(TX)\in H^{n}(X;\Z)$ is the \emph{Euler class} of X. We thus have $i^{*}(y) = \psi^{-1}(y) \smallsmile e(X)$ for all $y\in H^{*}(TX)$, and we obtain
$$\psi^{-1}\circ \ch(x) \smallsmile e(X) = i^{*}\bigl( \ch(x) \bigr) = \ch\bigl( i^{*}(x) \bigr).$$
Then, letting $D^{\bullet}$ be the elliptic complex over a compact manifold $X$ given by
$$
\xymatrix{
0 \ar[r] & \Gamma E^{0} \ar[r]^-{d_{0}} & \Gamma E^{1} \ar[r]^-{d_{1}} & \cdots
\ar[r]^-{d_{n-1}} & \Gamma E^{n} \ar[r] & 0\,,
}
$$
we see that $i^{*}\sigma(D^{\bullet}) = \sum_{i=0}^{n}(-1)^{i}E^{i}$, and so we have
$$\psi^{-1}\circ \ch\bigl( \sigma(D^{\bullet}) \bigr) \smallsmile e(X)
= \ch \left( \sum_{i=0}^{n}(-1)^{i}E^{i} \right).$$
We now incorporate this result into the statement of the Atiyah-Singer Index Theorem.

\begin{theorem}[Index Theorem C]
If $D^{\bullet}$ is an elliptic complex on a compact oriented $n$-dimensional manifold $X$ then
$$\chi(D^{\bullet}) = (-1)^{n(n+1)/2} \left( 
\ch\left( \sum_{i=0}^{n}(-1)^{i}E^{i} \right) \smallsmile e(X)^{-1} \smallsmile \td(TX\otimes\C) \right) [X],$$
provided that the element $\ch \bigl(\sum_{i=0}^{n}(-1)^{i}E^{i}\bigr)\smallsmile e(X)^{-1}\in H^{*}(X;\Q)$ is well defined.
\end{theorem}

We note that in this special case, the index of the elliptic complex depends only on the vector bundles involved and not on the action of the operators themselves. However, this version of the index theorem does not apply in general. In particular, we note that the Euler class $e(X)$ vanishes if $X$ admits a non-vanishing field of tangent vectors, and so $e(X)^{-1}$ is not defined.

\subsection{The Index on Odd-Dimensional Manifolds}

Let $D:\Gamma E \to \Gamma F$ be an elliptic partial differential operator over a compact, odd-dimensional manifold $X$. We have already noted that the symbol of an elliptic partial differential operator is a positively homogeneous map. In terms of local coordinates, we then have $\sigma(x,\lambda\xi) = \lambda^{m}\sigma(x,\xi)$ for all real $\lambda > 0$, where $m$ is the degree of the operator. However, since the symbol is locally a polynomial, it satisfies the even stronger condition that $\sigma(x,\lambda\xi) = \lambda^{m}\sigma(x,\xi)$ for \emph{any} $\lambda$, including negative values. In particular, we have $\sigma(x,-\xi) = \pm \sigma(x,\xi),$ where the sign is constant over the entire manifold. So, letting $\alpha : TX \to TX$ denote the fiberwise antipodal map  (bundle involution) $\alpha : (x,\xi) \mapsto (x,-\xi)$ given by multiplication by $-1$ on each fiber, we consider the induced map $\alpha^{*}: K(TX) \to K(TX)$. Representing the symbol class $\sigma(D)\in K(TX)$ by the complex
$$
\xymatrix{
0 \ar[r] & \pi^{*}E \ar[r]^{\sigma(x,\xi)} & \pi^{*}F \ar[r] & 0
}
$$ 
we see that the element $\alpha^{*}\sigma(P)$ is represented by the complex
$$
\xymatrix{
0 \ar[r] & \pi^{*}E \ar[r]^{\sigma(x,-\xi)} & \pi^{*}F \ar[r] & 0\,.
}
$$ 
We then obtain $\alpha^{*}\sigma(P) = \sigma(P),$ which follows immediately if $\sigma(x,-\xi) = + \sigma(x,\xi)$, while if $\sigma(x,-\xi) = -\sigma(x,\xi)$, then we can construct a homotopy between the two complexes by rotating halfway around the unit circle.

\begin{theorem}
The index of any elliptic partial differential operator on a compact odd-dimen\-sional manifold is zero.
\end{theorem}

\begin{proof}[Proof 1] Recalling that the Chern character induces a natural isomorphism from rational $K$-theory to rational cohomology, we have the following commutative diagram:
$$
\xymatrix{
K^{*}(TX)\otimes \Q \ar[r]^{\alpha^{*}} \ar[d]^{\ch}_{\cong} & K^{*}(TX)\otimes \Q \ar[d]^{\ch}_{\cong} \\
H_{c}^{*}(TX;\Q) \ar[r]^{\alpha^{*}} & H_{c}^{*}(TX;\Q)
}
$$
By the Thom isomorphism theorem of cohomology, we know that $H_{c}^{*}(TX;\Q)$ is a module over $H^{*}(X;\Q)$ generated by the Thom class $u$, where $u$ is the unique cohomology class that restricts to the preferred orientation class on each fiber of $TX$. Then, since $X$ is odd-dimensional, the antipodal map $\alpha : TX \to TX$ reverses the orientation on each fiber, and we obtain $\alpha^{*}u = -u$. It then follows that $\alpha^{*}x = -x$ for each $x\in H^{*}_{c}(TX;\Q)$, and so the corresponding map in $K$-theory is also given by $\alpha^{*}x = -x$ for each $x \in K^{*}(TX) \otimes \Q$. Since $\alpha^{*}\sigma(D) = \sigma(D)$, the image of $\sigma(D)$ in $K^{*}(TX)\otimes\Q$ is zero, and it follows that $\sigma(D)$ is an element of finite order in $K(TX)$. Since the topological index is a homomorphism from $K(TX)$ into the integers, we see that $\Index D = \tind \sigma(D) = 0.$
\end{proof}

\begin{proof}[Proof 2] By the cohomological form (A) of the index theorem, we have
$$\Index D = (-1)^{n}\bigl( \ch(\sigma(d)) \smallsmile \Im(X) \bigr) [TX],$$
where $n$ is the dimension of $X$. Applying the antipodal map $\alpha : TX \to TX$, we then obtain
\begin{equation*}\begin{split}
	\Index D &= (-1)^{n} \bigl( \ch (\alpha^{*}\sigma(D)) \smallsmile \Im(X) \bigr) \alpha_{*}[TX] \\
	&= (-1)^{n}\bigl( \ch(\sigma(D)) \smallsmile \Im(X) \bigr) (-1)^{n}[TX] \\
	&= -\Index D.
\end{split}\end{equation*}
We note that the index class $\Im(X) = \td(TX\otimes_{\R}\C)\in H^{*}(X,\Q)$ is invariant under the action of $\alpha^{*}$ because the complexification $\alpha:TX\otimes_{\R}\C \to TX \otimes_{\R}\C$ is homotopic to the identity by rotating halfway around the unit circle. It then follows that $\Index D = 0$.
\end{proof}

Note that in both of these proofs, the fundamental fact is that the antipodal map reverses the orientation of $TX$ in any odd-dimensional manifold. We also note that the first proof demonstrates that not every element of $K(TX)$ is the symbol of some elliptic partial differential operator. In particular, on the circle, we have $K(TS^{1}) \cong \tilde{K}(S^{2}) \cong \Z$, while every elliptic partial differential operator has symbol zero.

\subsection{The de\,Rham Complex}

We now take a second look at the (complexified) de\,Rham complex $D^{\bullet}$ given by
$$
\xymatrix{
0 \ar[r] & \Gamma \Lambda^{0}(TX\otimes \C) \ar[r]^-{d} &
\Gamma\Lambda^{1}(TX\otimes \C) \ar[r]^-{d} & \cdots \ar[r]^-{d} &
\Gamma\Lambda^{n}(TX\otimes \C) \ar[r] & 0
}
$$
By the de\,Rham theorem, we recall that the Euler characteristic of the complex is
$$\chi(D^{\bullet}) = \sum_{i}(-1)^{i}\dim H^{i}(D^{\bullet})
= \sum_{i}(-1)^{i}\dim H^{i}_{\mathrm{dR}}(X) = \chi(X).$$
First we consider the case where $X$ has even dimension $n=2l$. To compute the topological index, we note that by the elementary properties of the Euler class, we have
$$c_{n}(TX\otimes \C) = e\bigl( (TX\otimes\C)_{\R} \bigr) 
= (-1)^{l}e(TX\oplus TX) = (-1)^{l}e(X)^{2}.$$
where we must introduce the sign  $(-1)^{l}$ because of the difference in orientation between $(TX\otimes\C)_{\R}$ and $TX\oplus TX$.\footnote{Given an ordered basis $\{x_{1},\ldots,x_{n}\}$ for a real vector space $V$, we take
$\{ x_{1}.i x_{1},\ldots,x_{n},i x_{n}\}$ as our basis for $(V\otimes\C)_{\R}$, while we take $\{x_{1},\ldots,x_{n},x'_{1},\ldots,x'_{n}\}$ as our basis for $V\oplus V$. If $n=2l$, we then see that the permutation between these two bases then has sign $(-1)^{l}$.} We now apply a slight extension of the splitting principle which allows us to assume (in the even-dimensional case) that $TX \otimes \C$ decomposes as a direct sum of line bundles of the form $TX\otimes \C = (L_{1}\oplus \overline{L}_{1} ) \oplus \cdots \oplus (L_{l}\oplus \overline{L}_{l} ).$ Writing $c_{1}(L_{i}) = x_{i},$ we note that $c_{1}(\overline{L}_{i}) = -x_{i}$, and so we can calculate:
\begin{align*}
c_{n}(TX\otimes\C) &= \prod_{i}x_{i}\smallsmile (-x_{i})
= (-1)^{l}\prod_{i}x_{i}^{2}, \\
\Im(X) &= \td(TX\otimes \C) = \prod_{i}\frac{x_{i}}{1-e^{-x_{i}}}\frac{-x_{i}}{1-e^{x_{i}}}, \\
\ch \left( \sum_{i=0}^{n} (-1)^{i}\Lambda^{i}(TX\otimes\C) \right)
&= \ch \left( \prod_{i}(1-L_{i})(1+L_{i}) \right) \\
&= \prod_{i}(1-e^{x_{i}}) (1-e^{-x_{i}}) \\
&= c_{n}(TX\otimes \C) \smallsmile \td(TX\otimes \C)^{-1} \\
&= (-1)^{l}e(X)^{2} \smallsmile \td(TX\otimes \C)^{-1}.
\end{align*}
Noting that
$$\ch\left( \sum_{i=0}^{n}(-1)^{i}\Lambda^{i}(TX\otimes\C) \right)
\smallsmile e(X)^{-1} = (-1)^{l}e(X)\smallsmile \td(TX\otimes\C)^{-1}$$
is well defined even in the case where $e(X) = 0$, we can apply the cohomological form (C) of the Atiyah-Singer Index Theorem to obtain
$$\chi(X) = \chi(D^{\bullet}) = (-1)^{n(n+1)}(-1)^{l}e(X)[X] = e(X)[X].$$
In the odd-dimensional case, we note that the Euler class $e(X)$ vanishes, and by our above discussion of the index of differential operators on odd-dimensioanl manifolds, we see that $\chi(D^{\bullet}) = 0$ as well. Hence, for any manifold $X$, the Atiyah-Singer Index Theorem for the de\,Rham complex gives us $\xi(X) = e(X)[X].$

\subsection{The Dolbeault Complex}

We now examine the complex analog of the de\,Rham complex. Given a complex $n$-dimen\-sional manifold $X$, we let $T_{\C}X$ denote the complex tangent bundle of $X$ (i.e., the tangent bundle with respect to the complex structure on $X$), and we have the canonical isomorphism $TX\otimes \C \cong T_{\C}X \oplus \overline{T_{\C}X}$. We can then obtain the decomposition
$$\Lambda^{i}(TX\otimes\C) \cong \sum_{p+q=i}\Lambda^{p}(\overline{T_{\C}X})\otimes \Lambda^{q}(T_{\C}X) = \sum_{p+q=i}\Lambda^{p,q}(T_{\C}X),$$
where we define
$$\Lambda^{p,q}(T_{\C}X) = \Lambda^{p}(\overline{T_{\C}X}) \otimes \Lambda^{q}(T_{\C}X).$$
From this we see that the space $\Gamma \Lambda^{i}(TX\otimes\C)$ of complex differential forms on $X$ admits the decomposition
$$\Gamma \Lambda^{i}(TX\otimes\C) = \sum_{p+q=i}\Gamma\Lambda^{p,q}(T_{\C}X),$$
where $\Gamma\Lambda^{p,q}(T_{\C}X)$ is called the space of smooth differential forms of type $(p,q)$.

Taking local complex coordinates $z_{1},\ldots,z_{n}$ on $X$, where $z_{j} = x_{j} + i y_{j}$, the differential forms of type $(p,q)$ are given by
$$\omega = \sum_{\substack{j_{1}< \cdots < j_{p}\\k_{1}< \cdots< k_{q}}} a_{J,K}\,dz_{j_{1}}\cdots dz_{j_{p}}\,d\bar{z}_{k_{1}}\cdots d\bar{z}_{k_{q}},$$
where  $J = (j_{1},\ldots,j_{p})$ and $K = (k_{1},\ldots,k_{q})$ are multi-indices and $\bar{z}_{j} = x_{j} - i y_{j}$ denotes the complex conjuugate. The de\,Rham operator $d : \Gamma\Lambda^{i}(TX\otimes\C) \to \Gamma\Lambda^{i+1}(TX\otimes\C)$ then decomposes as the sum $d = \partial + \bar{\partial}$, where we define $\partial: \Gamma\Lambda^{p,q}(T_{\C}X) \to \Gamma\Lambda^{p+1,q}(T_{\C}X)$ and $\bar{\partial}: \Gamma\Lambda^{p,q}(T_{\C}X) \to \Gamma\Lambda^{p,q+1}(T_{\C}X)$ locally by
\begin{align*}
	\partial & : \sum_{J,K}a_{J,K}\,dz_{J}\,d\bar{z}_{K} \mapsto
	\sum_{j,J,K}\frac{\partial}{\partial z_{j}}\,a_{J,K}\,dz_{j}\,dz_{J}\,d\bar{z}_{k}, \\
	\bar{\partial} & : \sum_{J,K}a_{J,K}\,dz_{J}\,d\bar{z}_{K} \mapsto
	\sum_{k,J,K}\frac{\partial}{\partial \bar{z}_{k}}\,a_{J,K}\,d\bar{z}_{k}\,dz_{J}\,d\bar{z}_{k},
\end{align*}
where
$$\frac{\partial}{\partial z_{j}} = \frac{1}{2} \left( \frac{\partial}{\partial x_{j}} - i \frac{\partial}{\partial y_{j}} \right), \qquad
\frac{\partial}{\partial \bar{z}_{j}} = \frac{1}{2} \left( \frac{\partial}{\partial x_{j}} + i \frac{\partial}{\partial y_{j}} \right).$$
The \emph{Dolbeault complex} is then defined to be the elliptic complex
$$
\xymatrix{
	0 \ar[r] & \Gamma\Lambda^{0,0}(T_{\C}X) \ar[r]^-{\bar{\partial}}
	& \Gamma\Lambda^{0,1}(T_{\C}X) \ar[r]^-{\bar{\partial}}
	& \cdots \ar[r]^-{\bar{\partial}}
	& \Gamma\Lambda^{0,n}(T_{\C}X) \ar[r]
	& 0
}
$$
given by restricting the de\,Rham complex to differential forms of type $(0,q)$.

To compute the index of the Dolbeault complex, we note that since $X$ is a complex manifold, the Euler class is given by $e(X) = e(TX) = c_{n}(T_{\C}X)$.
It follows that we may apply form (C) of the index theorem. We recall that $\Lambda^{0,q}(T_{\C}X) = \Lambda^{q}(T_{\C}X)$ by definition. We can assume by applying the splitting principle that we have the decomposition $T_{\C}X = L_{1}\oplus \cdots \oplus L_{n}$, with $x_{i} = c_{1}(L_{i}),$ and we compute:
\begin{align*}
	\ch \left( \sum_{i=0}^{n}\Lambda^{i}(T_{\C}X) \right)
	&= \ch \left( \prod_{i}(1-L_{i}) \right) 
    = \prod_{i}(1-e^{x_{i}}), \\
    e(X) &= c_{n}(T_{\C}X) = \prod_{i}x_{i}, \\
    \Im(X) &= \td(TX\otimes \C) = \td(T_{\C}X) \smallsmile \td(\overline{T_{\C}X})
    = \td(T_{\C}X) \smallsmile \prod_{i}\frac{-x_{i}}{1-e^{x_{i}}} .
\end{align*}
Plugging these expressions into the index theorem, we obtain
\begin{align*}
	\chi(\bar{\partial}) &= (-1)^{2n(2n+1)/2}\Bigl(
	\ch \Bigl( \sum_{i=0}^{n}(-1)^{i}\Lambda^{i}(T_{\C}X) \Bigr)
	\smallsmile e(X)^{-1}\smallsmile \td(TX\otimes \C) \Bigr) [X] \\
	&= (-1)^{n(2n+1)}\left( \prod_{i}\frac{1-e^{x_{i}}}{x_{i}}
	\smallsmile \td(T_{\C}X) \smallsmile \prod_{i}\frac{-x_{i}}{1-e^{x_{i}}} \right)
	[X] \\
	&= (-1)^{n}\bigl( (-1)^{n}\td(T_{\C}X)\bigr) [X] = \td(T_{\C}X) [X].
\end{align*}
By a generalization of this construction, we can obtain the celebrated Hirzebruch-Riemann-Roch theorem. For a complete discussion of this, see \cite[\S 4]{AS68III}.

\begin{example}
Viewing the sphere $S^{2}$ as a complex manifold of dimension 1 (i.e., the Riemann sphere), the Dolbeault complex reduces to the Dolbeault operator
$$\bar{\partial} : \Gamma \Lambda^{0}(T_{\C}S^{2}) \to \Gamma\Lambda^{1}(T_{\C}S^{2})$$
over $S^{2}$ given by $\bar{\partial}f = ( \partial / \partial \bar{z}) \, d\bar{z}.$ By the above discussion, we then have
$$\Index \bar{\partial} = \td(T_{\C}S^{2})[X].$$
Now, recalling that the Todd class is given by
$$\td(E) = 1 + \frac{1}{2}\,c_{1}(E) + \text{(higher order terms),}$$
we obtain
\begin{align*}
	\Index\bar{\partial} &= \td(T_{\C}S^{2})[X]
	= \bigl( (1+ \frac{1}{2}\,c_{1}(T_{\C}S^{2}) \bigr) [X] \\
	&= \frac{1}{2} \bigl( c_{1}(T_{\C}S^{2}) \bigr) [X]
	= \frac{1}{2}\, e(X) [X],
\end{align*}
and since we just proved that $e(X)[X] = \chi(X)$, we see that
$$\Index \bar{\partial} = \frac{1}{2}\,\chi(S^{2}) = 1.$$
More generally, letting $X$ be a Riemann surface of genus $g$, we have $\chi(X) = 2 - 2g$, and so we obtain $\Index \bar{\partial}_{X} = 1 - g.$
\end{example}

\section{Proof of the Index Theorem}
\label{sec:3}

\subsection{Axioms for the Topological Index}

It follows immediately from our construction of the topological index that it satisfies the following two elementary axioms:
\begin{description}
\item[Axiom A] If $X$ is a point, then $\tind_{\pt}:\Z \to \Z$ is the identity map.
\item[Axiom B] $\tind$ commutes with the homomorphisms $i_{!}$.
\end{description}

In fact, we shall see that these two axioms uniquely characterize the topological index.

\begin{definition}
A collection of homomorphisms $\ind_{X} : K(TX) \to \Z$ for all compact manifolds $X$ is called an \emph{index function} if it is functorial with respect to diffeomorphisms. In other words, if $f : X\to Y$ is a diffeomorphism, then $\ind_{X}f^{*}y = \ind_{Y}y$ for all $y\in K(TY)$, where $f^{*}:K(TY) \to K(TX)$ is the map induced by the extension of $f$ to the tangent bundles.
\end{definition}

\begin{proposition}
If $\ind$ is an index function satisfying Axioms A and B above, then
$$\ind = \tind.$$
\end{proposition}

\begin{proof}
For any compact manifold $X$, we take an embedding $i : X \to \R^{m}$, and we let $j : \pt \to \R^{m}$ be the inclusion of the origin. By Axiom A above, we see that the diagram
$$
\xymatrix{
K(TX) \ar[r]^{i_{!}} \ar[dr]_{\ind} & K(T\R^{m}) \ar[r]^{j_{!}^{-1}} & K(T\pt) \ar[dl]^{\ind} \\
& \Z
}$$
is commutative, and so we obtain $\ind_{X} = \ind_{\pt}\circ j_{!}^{-1} \circ i_{!} = j_{!}^{-1} \circ i_{!} = \tind_{X}.$
\end{proof}

The direction of our proof of the Atiyah-Singer Index Theorem is now clear. We would like to construct an analytical index given at the symbolic level by $\aind = \Index \circ \sigma^{-1}$, mapping each element $x\in K(TX)$ to the index of some elliptic operator with symbol $x$. The index theorem then reduces to showing that the topological index and the analytical index coincide, and so we need only verify Axioms A and B for this analytical index. In this case, Axiom B tells us that given an elliptic operator $D$ on $X$ and an inclusion map $i : X \to  Y$, we can construct (symbolically) an elliptic operator $i_{!}D$ on $Y$ which has the same index as $D$. Then, since any manifold can be embedded in $\R^{m}$, we can reduce the index problem to that of elliptic operators over the $m$-sphere $S^{m} = (\R^{m})^{+}$, which is much more easily solved. Unfortunately, as we saw in our discussion of odd-dimensional manifolds in \S\ref{sec:2}, we are not guaranteed that every element of $K(TX)$ is the symbol of some elliptic partial differential operator. To fix this, we need to extend our discussion to the class of pseudo-differential operators.

\subsection{Pseudo-Differential Operators}

Let $f$ be a real-valued function on $\R^{n}$ with compact support. We define the Fourier transform of $f$ to be the function $\hat{f}$ given by the integral
$$\hat{f}(\xi) = \frac{1}{(2\pi)^{n}} \int_{\R^{n}}f(x)\,e^{-i\langle x,\xi\rangle}dx.$$
The Fourier inversion formula lets us write $f$ in terms of $\hat{f}$ as the integral
$$f(x) = \int_{\R^{n}}\hat{f}(\xi) \, e^{i\langle x,\xi\rangle}.$$
Taking the partial derivative $D_{j} = -i\,\partial/\partial x_{j}$ of $f$, we obtain
\begin{align*}
	-i\,\frac{\partial}{\partial x_{j}}\,f(x)
	&= \int_{\R^{n}}-i\,\frac{\partial}{\partial x_{j}}\,\hat{f}(\xi)\,e^{i\langle x,\xi\rangle}d\xi \\
	&= \int_{\R^{n}} \xi_{j}\hat{f}(\xi)\, e^{i\langle x,\xi\rangle}d\xi.
\end{align*}
Hence, the Fourier transform ``converts differentiation into multiplication''. In particular, if $P$ is a partial differential operator on $\R^{n}$, then we have
$$Pf(x) = \int_{\R^{n}}p(x,\xi)\hat{f}(\xi)\,e^{i\langle x,\xi\rangle}d\xi,$$
where $p(x,\xi)$ is a polynomial in $\xi = (\xi_{1},\ldots,\xi_{n})$ whose coefficients are smooth real-valued functions of $x$. To make clear the connection between $P$ and $p(x,\xi)$, we can write $P = p(x,D)$, where $D = (-i\,\partial/\partial x_{1}, \ldots, -i\,\partial/\partial x_{n} ).$ Note that we can ``extract'' the polynomial $p(x,\xi)$ from the operator  $P$ by taking the commutator
$$p(x,\xi) = e^{-i\langle x,\xi \rangle} P e^{i\langle x,\xi\rangle}.$$
The leap to pseudo-differential operators comes when we consider functions $p(x,\xi)$ which are not necessarily polynomials. All we require are suitable growth conditions on $p(x,\xi)$ and the ability to define the symbol $\sigma(x,\xi).$

\begin{definition}
A linear operator $P$ from smooth functions on $\R^{n}$ with compact support to smooth functions on $\R^{n}$ is called a \emph{pseudo-differential operator} of order $m$ if it is given by $P = p(x,D)$, where $p(x,\xi)$ is a smooth function satisfying the growth conditions\footnote{These growth conditions allow us to differentiate under the integral sign.}
$$ \bigl| D^{\beta}_{x}\,D^{\alpha}_{\xi}\, p(x,\xi) \bigr| 
\leq C_{\alpha,\beta}\bigl( 1 + |\xi| \bigr) ^{m-|\alpha|},$$
with a positive constant $C_{\alpha,\beta}$ for each multi-index $\alpha,\beta$. Here we write $|\alpha| = \alpha_{1} + \cdots + \alpha_{n}$ and $D^{\alpha}_{\xi} = (\partial/\partial\xi_{1})^{\alpha_{1}}\cdots (\partial/\partial\xi_{n})^{\alpha_{n}}$.  We also require that for all $x\in \R^{n}$ and $\xi \neq 0$, the limit
$$\sigma_{m}(x,\xi) = \lim_{\lambda\to \infty}\frac{p(x,\lambda\xi)}{\lambda^{m}}$$
exists, and $\sigma_{m}(x,\xi)$ is called the $m$-th order \emph{symbol} of $P$.
\end{definition}

If $P$ is a partial differential operator, then the corresponding polynomial $p(x,\xi)$ satisfies the growth conditions, and $\sigma_{m}(x,\xi)$ is equal to the symbol we defined by removing all but the highest order terms of $p(x,\xi)$. Also note that the order $m$ need not be an integer. We can consider operators with $m = \frac{1}{2}$ such as $P = p(x,D)$ where $p(x,\xi) = a_{1}(x) \sqrt{\xi_{1}} + \cdots a_{n}(x) \sqrt{\xi_{n}},$ and we can also consider operators with $m \leq 0$. In any case, the symbol $\sigma_{m}(x,\xi)$ of $P$ has the property that
$$\sigma_{m}(x,\lambda\xi) = \lambda^{m}\sigma_{m}(x,\lambda\xi) \text{ for all real } \lambda > 0,$$
and so we say that $\sigma_{m}(x,\xi)$ is \emph{positively homogeneous of degree $m$.}

One problem that we encounter with this definition is that the function $e^{i\langle x,\xi \rangle}$ does not have compact support, and so we cannot take the commutator given above. To remedy this, we localize our definition by requiring only that for each function $f$ with compact support there exists a smooth function $p_{f}(x,\xi)$ satisfying the growth conditions given above such that
$$P(f\cdot u) = p_{f}(x,D)\,u$$
for any function $u$ with arbitrary support. In this case, we can now take the commutator
$$p_{f}(x,\xi) = e^{-\langle x,\xi \rangle} P\bigl( f \cdot e^{i\langle x,\xi\rangle}\bigr).$$
Then, we can define the symbol of $P$ by $\sigma(x,\xi) = \sigma_{f}(x,\xi)$, where $f$ is a function with compact support that is equal to $1$ in some neighborhood of $x$. It turns out that this definition is independent of the choice of the function $f$.

Our definition easily generalizes to the case of pseudo-differential operators on vector-valued functions $f: \R^{m}\to \R^{n}$, in which case the polynomials $p_{f}(x,\xi)$ and the symbol $\sigma(x,\xi)$ are matrix-valued functions of $x$ and $\xi$. We say that the pseudo-differential operator $P$ is \emph{elliptic} if its symbol $\sigma(x,\xi)$ is invertible for all $x$ and $\xi\neq 0$. We are now ready to extend the concept of pseudo-differential operators to manifolds.

\begin{definition}
Let $E$ and $F$ be complex vector bundles over a manifold $X$. A linear operator\footnote{Here we use the notation $\Gamma_{c}$ to denote the space of smooth sections with \emph{compact support}.} $P: \Gamma_{c}E \to \Gamma F$ is called a \emph{pseudo-differential operator} of degree $m$ on $X$ if given a covering $\{U_{u}\}$ of $X$ by coordinate patches over which $E$, $F$, and $TX$ are trivial, the operators $P_{i}$ obtained by restricting $P$ to functions with compact support in $U_{i}$ are pseudo-differential operators of degree $m$ on $\R^{n}$ (where $n$ is the dimension of $X$).
\end{definition}

We can then define a global symbol $\sigma(P) : \pi^{*}E \to \pi^{*}F$, and we say that a pseudo-differential operator $P$ is \emph{elliptic} if its symbol $\sigma(P)$ is invertible outside the zero section. As before, we see that if $X$ is compact, the symbol gives us an element $\sigma(P)\in K(TX)$. In fact, we can obtain any element of $K(TX)$ in this manner. Before demonstrating this, we will first work out an example.

\begin{example}
In \S\ref{sec:2}, we showed that any elliptic partial differential operator on an odd dimensional manifold has index $0$. In particular, any such operator on the circle has index zero. We will now construct an elliptic pseudo-differential operator on the circle with index $-1$. consider the linear operator $P$ on complex-valued functions over the circle $S^{1} = \R / 2\pi\Z$ given in terms of the Fourier series expansion by
$$	P(e^{ik\theta}) = \begin{cases}
e^{ik\theta} & \text{ for $k \geq 0$,} \\
0 & \text{ for $k < 0$.}
\end{cases}$$
So $P$ essentially kills the negative portion of the Fourier series expansion. Given any complex-valued function $f$ with compact support in the interval $0 < \theta < 2 \pi$, we obtain
\begin{align*}
	p_{f}(\theta,\xi) &= e^{-i\theta\xi}P(f\cdot e^{i\theta\xi})
	= e^{-i\theta\xi} P \left( \sum_{-\infty}^{\infty}\hat{f}(k) \,e^{ik\theta}\cdot e^{i\theta\xi}\right) \\
	&= e^{-i\theta\xi} P \left( \sum_{-\infty}^{\infty}\hat{f}(k-\xi)\,e^{ik\theta}\right)
	= \sigma_{0}^{\infty}\hat{f}(k-\xi)\,e^{i\theta(k-\xi)},
\end{align*}
and taking the limit we see that
$$p_{f}(\theta,\xi) \to \begin{cases}
f(\theta) & \text{ as $\xi \to +\infty$,} \\
0 & \text{ as $\xi \to -\infty$.}
\end{cases}$$
It follows that $p_{f}(\theta,\xi)$ satisfies the growth conditions for a pseudo-differential operator of degree zero, where the global $0$-th order symbol is given by
$$\sigma_{p}(\theta,\xi) = \begin{cases}
1 & \text{ for $\xi > 0$,} \\
0 & \text{ for $\xi < 0.$}
\end{cases}$$
Although this operator is not elliptic, the operator $A = e^{i\theta}P + (1-P)$ is also pseudo-differential of degree zero, and since its symbol is
$$\sigma_{A}(\theta,\xi) = \begin{cases}
e^{i\theta} & \text{ for $\xi > 0$,} \\
1 & \text{ for $\xi < 0$.}
\end{cases}$$
we see that $A$ is elliptic. Furthermore, in terms of the Fourier series expansion, we obtain
$$A(e^{ik\theta}) = \begin{cases}
	e^{i(k+1)\theta} & \text{ for $k\geq 0$,} \\
	e^{ik\theta} & \text{ for $k < 0$.}
\end{cases}$$
The operator $A$ thus shifts the positive portion of the Fourier series by one term, and we see that $\Ker A$ is empty while $\Coker A$ can be identified with the space of constant functions. It follows that $\Index A = -1$. For a generalization of this example in the context of a proof of the Bott periodicity theorem, see the appendix.
\end{example}

\subsection{Construction of the Analytical Index}

\begin{claim}
Given an element $a\in K(TX)$, where $X$ is a compact manifold, there exists an elliptic pseudo-differential operator $P$ of arbitrary order $m$ on $X$ such that $\sigma(P) = a$.
\end{claim}

\begin{proof}
We know that $a$ can be represented by a complex of length one
$$\xymatrix{
0 \ar[r] & \pi^{*}E \ar[r]^{\alpha_{m}} & \pi^{*}F \ar[r] & 0
},$$
where $\alpha_{m}$ is homogeneous of degree $m$. Choosing local coordinates that trivialize $E$, $F$, and $TX$, we take $p(x,\xi) = \varphi(\xi)\alpha_{m}(x,\xi),$ where $\varphi(\xi)$ is a smooth function that is zero in a neighborhood of the origin and is $1$ elsewhere. This correction term is necessary since for $m \leq 0$ we observe that $\alpha_{m}$ may be discontinuous at the zero section. It is then easy to see that $P = p(x,D)$ is a pseudo-differential operator of order $m$ on $X$, and that its symbol is precisely $\alpha_{m}$. Hence, $P$ is an elliptic operator such that $\sigma(P) = a$.
\end{proof}

Now that we have shown that the symbol map is surjective, we would like to show that two elliptic operators with the same symbol have the same index (we could then say that the symbol map is ``pseudo-injective''). Before we show this, we must be sure that the index of an elliptic pseudo-differential operator is actually defined. Fortunately, if $P$ is an elliptic pseudo-differential operator, we know from analysis that $\Ker P$ and $\Coker P$ are both finite-dimensional\footnote{This is not at all an immediate conclusion. See \cite[\S 5]{AS68I} for a discussion involving Sobolev spaces. If it still bothers you, you should discuss this problem with an analyst.} and so we may define the index to be
$$\Index P = \dim \Ker P - \dim\Coker P.$$
We also know that the index of an elliptic pseudo-differential operator of order $m$ depends on the homotopy class of the symbol in the space of symbols with the same order. We are now ready to prove:

\begin{claim}
If $P$ and $P'$ are two elliptic pseudo-differential operators on a compact manifold $X$ such that $\sigma(P)= \sigma(P')\in  K(TX)$, then $\Index P = \Index P'$.
\end{claim}

\begin{proof}
Suppose that $\sigma(P)$ and $\sigma(P')$ are given by the complexes
$$\xymatrix{
0 \ar[r] & \pi^{*}E \ar[r]^{\sigma_{m}} & \pi^{*}F \ar[r] & 0}
\qquad
\xymatrix{
0 \ar[r] & \pi^{*}E' \ar[r]^{\sigma'_{n}} & \pi^{*}F' \ar[r] & 0},
$$
where $\sigma_{m}$ and $\sigma'_{n}$ are the symbols of $P$ and $P'$ respectively. For now, we suppose that $P$ and $P'$ are both of order $m$. Since $\sigma(P) = \sigma(P')$, we see that the two complexes are either homotopic or else they differ by a complex with empty support. If they are homotopic, then $\sigma_{m}$ and $\sigma'_{m}$ are homotopic via homogeneous homomorphisms of degree $m$, and it follows that $\Index P = \Index P'$. If they differ by a complex of the form
$$\xymatrix{
0 \ar[r] & G \ar[r]^{\beta} & G \ar[r] & 0
},$$
where $\beta$ is an isomorphism, then since $\beta$ is also homogeneous of degree $m$, we have $\dim G = 0$ unless $m=0$. When $m=0$, we see that $\beta$ is constant in $\xi$, and so we see that $P$ and $P'$ differ by the operator $Qu = \beta(x,0)u$ with $\Index Q = 0$. It again follows that $\Index P = \Index P'$.

Now, we consider the case where $P$ is of order $m$ and $P'$ is of order $n$ with $m\neq n$. By the argument above, we may assume that $E = E'$ and $F=F'$ and that $\sigma_{m} = \sigma'_{n}$ on $S(X)$. Then, if we construct the map $\rho = \sigma'_{n}\sigma_{m}^{-1}$, we see that $\rho$ is homogeneous of degree $n-m$, and that $\rho$ is the identity map on $S(X)$. The map $\rho$ is thus self-adjoint, and it follows that $P$ and $P'$ differ by an associated self-adjoint operator $R$. Since $R$ is self-adjoint, we see that $\Index R = \dim \Ker R - \dim \Ker R^{*} = 0$, and so we obtain $\Index P = \Index P'$.
\end{proof}

\begin{definition}
The \emph{analytical index} is the map $\aind:K(TX) \to \Z$ given by
$$\aind : x \mapsto \Index P, \text{ where } \sigma(P) = x.$$
\end{definition}

\subsection{Verification of the Axioms}

From its construction, it is easy to see that $\aind$ is functorial with respect to diffeomorphisms, and so $\aind$ is an index function. We would then like to show that the analyticial index satisfies the following two axioms:
\begin{description}
\item[Axiom A] If $X$ is a point, then $\aind {\pt}:\Z \to \Z$ is the identity map.
\item[Axiom B] $\aind$ commutes with the homomorphisms $i_{!}$.
\end{description}
Since these two axioms uniquely characterize the topological index, $\tind$, we would then obtain $\aind=\tind$, and the Atiyah-Singer Index Theorem follows immediately. To verify Axiom A, we note that if $X$ is a point, then an elliptic operator on $X$ is simply a linear transformation $P : V \to W$ between two complex vector spaces $V$ and $W$, and
$$\sigma(P) = \dim E - \dim F \in K(TX)\cong \Z.$$
Furthermore, we note that
$$\dim \Ker P+ \dim \im P = \dim E, \qquad
  \dim F - \dim \im P = \dim \Coker P.$$
We thus have
$$\Index P = \dim\Ker P - \dim\Coker P = \sigma(P),$$
and Axiom A follows immediately.

The verification of Axiom B is not quite so simple. We recall that for an inclusion $i : X \to Y$ of a compact submanifold $X$ in $Y$, we defined $i_{!}$ to be the composition
$$ i_{!} : K(TX) \xrightarrow{\varphi}K(TN) \xrightarrow{h_{*}}K(TY),$$
where $N$ is a tubular neighborhood of $X$ in $Y$, $\varphi$ is the Thom isomorphism, and $h_{*}$ is the natural extension homomorphism induced by the open inclusion $h : TN \to TY$. To verify that $\aind$ commutes with $i_{!}$, we must show that it commutes with $\varphi$ and $h_{*}$. First we consider the natural extension homomorphism $h_{*}$.

\begin{proposition}[Excision]
If $h : U \to X$ is the inclusion of an open set $U$ in $X$, then
$$\aind \circ h_{*}(x) = h_{*}\circ \aind(x)$$
where $h_{*}:K(TU) \to K(TX)$ is the natural extension homomorphism.
\end{proposition}

\begin{proof} Any element $a\in K(TU)$ can be represented by a homogeneous complex
$$\xymatrix{
0 \ar[r] & \pi^{*}E \ar[r]^{\alpha} & \pi^{*}F \ar[r] & 0
}$$
over $TU$. Since $U$ is not necessarily compact, we must take $\alpha$ to be homogeneous of degree \emph{zero}, and we also require that $\alpha(x,\xi) = \text{Identity}$ for all $x$ outside some compact set $C\subset U$ (see \cite[pp.~493--3]{AS68I} for a discussion of this non-compact case). We can then trivially extend $E$ and $F$ to bundles $E'$ and $F'$ over $X$, and we consider the complex
$$\xymatrix{
0 \ar[r] & \pi^{*}E' \ar[r]^{\alpha} & \pi^{*}F' \ar[r] & 0
}$$
over $X$ representing the element $h_{*}a\in K(TX)$, where
$$\alpha' = \begin{cases}
\alpha & \text{ on $TU$,} \\
\text{Identity} & \text{ on $TX - TU$.}
\end{cases}$$
Letting $P$ and $P'$ be the corresponding operators over $U$ and $X$ respectively, then we have $\sigma(P') = h_{*}\sigma(P)$ by construction.
If $f$ has support in $U$ and $Pf = 0$, it follows that $f$ must have support in $U$ and so $Pf = 0$. Hence $\Ker P = \Ker P'$. The same is true for the adjoints $P^{*}$ and $P'^{*}$, and so we obtain $\Index P = \Index P'$.
\end{proof}

It only remains to show that the analytical index commutes with the Thom isomorphism $\varphi : K(X) \to K(V)$, where $V$ is a complex vector bundle over $X$. We recall that $\varphi(x) = \pi^{*}x \cdot \lambda_{V}$, where $\pi : V\to X$ is the vector bundle projection and $\lambda_{V} = \varphi(1)$ is the Thom class given by the exterior complex $\Lambda(V)$. We would like to show that
$$\aind ( \pi^{*}x\cdot \lambda_{V}) = \aind(x) \cdot \aind(\lambda_{n}) = \aind(x).$$
In order to do this, it is necessary to establish a product formula for $\aind$ and then verify that $\aind(\lambda_{n}) = 1$ under the appropriate circumstances. We begin with the product formula. This problem is considerably simplified if we consider trivial bundles which can be expressed as the product $V = X \times \R^{n}$. For general vector bundles we must instead consider ``twisted'' products.

Suppose that we have a principal fiber bundle $P\to X$ with an associated compact Lie group $H$, where $H$ acts freely on $P$ on the right and $X = P/H$. Then, given a locally compact space $F$ on  which $H$ acts on the left, we form the fiber bundle $Y$ over $X$ given by
$$Y = P \times_{H} F = (P\times F)\, /\, (p,f)\sim(ph^{-1},hf) \,\forall h\in H.$$
The action of $H$ on $F$ induces a $H$-action on $TF$, and so we can form the vector bundle $P\times_{H}TF$ over $Y$. Choosing a metric, we then note that we have the decomposition
$TY = (P\times_{H}TF) \oplus \pi^{*}TX$, which gives us a product map $K(TX) \otimes K(P\times_{H}TF) \to K(TY).$ Considering the homomorphism
$$K_{H}(TF) \to K_{H}(P\times TF) \cong K(P\times_{H}TF),$$
where $K_{H}(\cdot)$ denotes the \emph{equivariant $K$-theory} functor discussed in \S\ref{sec:4}, we obtain a multiplication map
$$K(TX) \otimes K_{H}(TF) \to K(TY).$$
For our purposes, we are interested in the case where $H = O(n)$ and $F = \R^{n}$. Then, any vector bundle $Y$ over $X$ can be written in the form $Y = P \times_{O(n)}\R^{n}$ for an appropriate principal fiber bundle $P\to X$, and we obtain the multiplication map
$$K(TX) \otimes K_{O(n)}(T\R^{n}) \to K(TY).$$
We are now ready to present our product formula for the analytical index.

\begin{proposition}[Multiplicative Axiom] If $a\in K(TX)$ and $b\in K_{H}(TF)$, then
$$\aind^{Y} a\cdot b = \aind^{H}a \cdot \aind_{H}^{F}b,$$
provided that $\aind^{F}_{H}b\in \Z$ is a multiple of the trivial representation $1\in R(H)$, where we define $\aind_{H}b = \Index_{H}P = [\Ker P] - [\Coker P]\in R(H)$ for $\sigma(P) = b$ (see \S\ref{sec:4}).
\end{proposition}

\begin{proof}
We can represent $a$ by a homogeneous complex of degree $1$ over $TX$
$$\xymatrix{
0 \ar[r] & \pi^{*}E^{0} \ar[r]^{\alpha} & \pi^{*}E^{1} \ar[r] & 0,
}$$
and we construct a pseudo-differential elliptic operator $A$ over $X$ with symbol $\alpha$. Then, by taking a partition of unity subordinate to an open covering $\{U_{i}\}$ which trivializes $Y$, we can lift $A$ to a global elliptic operator $\tilde{A}$ on $Y$.\footnote{If $P : \Gamma E \to \Gamma F$ is an operator on $X$ and $G$ is a vector bundle over $Y$, then we define the lifting $\tilde{P} : \Gamma ( E \tensorhat G) \to \Gamma(F \tensorhat G)$ of $P$ to $X\times Y$ by $\tilde{P}\bigl( u(x) \otimes v(y) \bigr) = P \bigl( u(x) \bigr) \otimes v(y).$} Similarly, we can represent $b$ by a homogeneous complex of degree $1$ over $TF$
$$\xymatrix{
0 \ar[r] & \pi^{*}G^{0} \ar[r]^{\beta} & \pi^{*}G^{1} \ar[r] & 0
}.$$
Letting $B$ be an elliptic pseudo-differential operator over $F$ with symbol $\beta$ (and furthermore requiring that $B$ commute with the actions of $H$ on $G^{0}$ and $G^{1}$), we can lift $B$ to a global elliptic operator $\tilde{B}$ on $Y$.

We now construct the operator $D$ on $Y$ given by
$$ D = \begin{pmatrix}
\tilde{A} & -\tilde{B}^{*} \\ \tilde{B} & \phantom{-}\tilde{A}^{*}
\end{pmatrix}, \text{ with }
\sigma(D) = \begin{pmatrix}
\alpha \tensorhat 1 & -1\tensorhat \beta^{*} \\ 1\tensorhat \beta  & \alpha^{*}\tensorhat 1
\end{pmatrix}.$$
Recalling the expression for the product of two complexes of length one, we note that $\sigma(D)$ corresponds to the element $ab\in K(TY)$ represented by the complex
$$\xymatrix{
0 \ar[r] & E^{0}\tensorhat G^{0}\oplus E^{1}\tensorhat G^{1} \ar[r]^{\sigma(D)}
& E^{1}\tensorhat G^{0} \oplus E^{0}\tensorhat G^{1} \ar[r] & 0
},$$
and so we must calculate the index of $D$. To do this, we consider the following two diagonalized operators (noting, of course, that $\tilde{A}$ and $\tilde{B}$ commute)
\begin{align*}
	D^{*}D &= \begin{pmatrix}
	\tilde{A}^{*}\tilde{A} + \tilde{B}^{*}\tilde{B} & 0 \\
	0 & \tilde{A}\tilde{A}^{*} + \tilde{B}\tilde{B}^{*}
	\end{pmatrix} = \begin{pmatrix}
	P_{0} & 0 \\ 0  & Q_{0}
	\end{pmatrix}, \\
	DD^{*} &= \begin{pmatrix}
	\tilde{A}\tilde{A}^{*} + \tilde{B}^{*}\tilde{B} & 0 \\
	0 & \tilde{A}^{*}\tilde{A} + \tilde{B}\tilde{B}^{*}
	\end{pmatrix} = \begin{pmatrix}
	P_{1} & 0 \\ 0  & Q_{1}
	\end{pmatrix}.
\end{align*}
where
\begin{align*}
\Ker D &= \Ker D^{*}D = \Ker P_{0}\oplus \Ker Q_{0}, \\
\Coker D &= \Ker D^{*} = \Ker DD^{*} = \Ker P_{1}\oplus \Ker Q_{1}.
\end{align*}
We thus have
$$\Index D = ( \dim \Ker P_{0} - \dim \Ker P_{1} )
+ ( \dim \Ker Q_{0} - \dim\Ker Q_{1} ).$$
Considering the operator $P_{0} = \tilde{A}^{*}\tilde{A} + \tilde{B}^{*}\tilde{B}$, we note that
$$\langle P_{0}u,u\rangle = \langle \tilde{A}u,\tilde{A}u\rangle
+ \langle \tilde{B}u, \tilde{B}u \rangle,$$
and so we have $\Ker P_{0} = \Ker \tilde{A} \cap \Ker \tilde{B}$. Since $\tilde{B}$ extends $B$ to the fibers of $Y$, we see that $\Ker \tilde{B}$ is the space of smooth sections of the vector bundle $K_{B} = P \times_{H}\Ker B$
over $X$. Then $\tilde{A}$ induces an operator $C$ on sections of $K_{B}$ with $\sigma(C) = \alpha \otimes \mathrm{Id}(K_{B})$, and it follows that $\sigma(C) = a[K_{B}]\in K(TX)$, where $[K_{B}]$ is the class of $K_{B}$ in $K(X)$. Replacing $\tilde{A}$ by $\tilde{A}^{*}$, we obtain the analogous result for $P_{1}$ and $C^{*}$, which gives us
\begin{align*} 
\dim\Ker P_{0} - \dim\Ker P_{1} &= \dim \Ker C - \dim \Ker C^{*} \\
&= \Index C = \aind^{X}(a[K_{B}]).
\end{align*}
Similarly, taking $L_{B} = P \times_{H}\Coker B$, we obtain
$$\dim \Ker Q_{0} - \dim \Ker Q_{1} = \aind^{X}(a[L_{B}]).$$
Combining these two results, we have
$$\Index D = \aind^{X} \bigl( a \cdot ([K_{B}] - [L_{B}]) \bigr) .$$
If in addition $\aind_{H}^{F}b = [\Ker B] - [\Coker B] \in R(H)$ is an integer (i.e., it is a multiple of the trivial representation $\C$), then $[K_{B}] - [L_{B}] = (\aind_{H}^{F}b)\cdot \mathbf{1}$ is a multiple of the trivial bundle $\mathbf{1}$, and since $\aind^{X}$ is a homomorphism we obtain our desired product formula
$\Index^{Y}D = (\aind^{X}a ) \cdot  (a\ind^{F}_{H}b).$
\end{proof}

\begin{proposition}[Normalization Axiom]
If $j : \pt \to \R^{n}$ is the inclusion of the origin, so that we have the induced homomorphism $j_{!}: R\bigl(O(n)\bigr) \to K_{O(n)}(T\R^{n})$, then
$$\aind_{O(n)}j_{!}(1) = \aind_{O(n)}\lambda_{n} = 1,$$
where the Bott class $\lambda_{n} = \lambda_{\C^{n}}$ is the $K$-theory Thom class of $\C^{n}$.
\end{proposition}

\begin{proof}
Considering first the case where $n=2$, we recall that the de\,Rham symbol $\rho_{2} = \rho_{S^{2}}$ is given by the complexification of the de\,Rham complex for the $2$-sphere $S^{2}$. Since the bundles involved are orientable, we may let the structure group be $SO(2)$. Considering the action of $SO(2)$ on $S^{2} = \R^{2}\cup \infty$, we see that $\rho_{2}\in K_{SO(2)}(TS^{2})$, and by a suitable deformation we obtain (see \cite[\S 3.2]{AS68I})
$$\rho_{2} = h(\lambda_{2}) + f^{*}\circ h(\lambda_{2}),$$
where $h$ is the extension homomorphism $K(T\R^{2}) \to K(TS^{2})$ and $f : S^{2}\to S^{2}$ is the reflection exchanging $0$ and $\infty$. Then since $\aind$ is functorial and commutes with the natural extension homomorphism, we have $\aind_{SO(2)} \rho_{2} = 2\aind_{SO(2)}\lambda_{2}.$ Noting that $\aind_{SO(2)} \lambda_{2} = \chi(S^{2}) = 2$, we then have $\aind_{SO(2)}\lambda_{2} = 1.$

In the case where $n=1$, we recall the elliptic operator $A$ on the circle $S^{1}$ from our discussion of pseudo-differential operators. By an appropriate homotopy, we then obtain  $-h(\lambda_{1}) = \sigma(A)$, where $h$ is the natural extension homomorphism $K(T\R^{1}) \to K(TS^{1})$. Recalling that $\Index A = -1$, we see that $\aind_{O(1)}\lambda_{1} = 1$. This construction is carried out in detail as an example of our discussion of clutching functions in the appendix.

Now, taking a decomposition of $\R^{n}$ into factors of $\R^{1}$ and $\R^{2}$, we see by the product formula that $\aind_{G}\lambda_{n}=1$ for all subgroups $G$ of $O(n)$ given by $G_{1}\times \cdots \times G_{n}$, where each $G_{i}$ is either $O(1)$ or $SO(2)$. Since a representation of $O(n)$ is completely determined by its restriction to these subgroups, we see that $\aind_{O(n)}\lambda_{n} = 1.$
\end{proof}

So, if we let $\varphi : K(TX) \to K(TY)$ be the Thom isomorphism for a complex vector bundle $TY$ over $TX$, we can combine the above two propositions to obtain
\begin{align*}
\aind^{Y}\circ \varphi(a) &= \aind^{Y}(\pi^{*}a \cdot \lambda_{TY}) \\
&= (\aind^{X}a ) \cdot \bigl(\aind^{R^{n}}_{O(n)}\lambda_{n}\bigr)
= \aind^{X}a.
\end{align*}
Now that we have shown that the analytical index commutes with the Thom isomorphism, our proof of the Atiyah-Singer Index Theorem is complete.
 
\section{Equivariant $K$-Theory and Homogeneous Spaces}
\label{sec:4}

\subsection{$G$-Vector bundles and $K_{G}(X)$}

We now consider the generalization of $K$-theory where we take into account the action of a compact Lie group $G$. Recall that a \emph{Lie group} is a group with the structure of a smooth manifold where the group multiplication and inverse maps are smooth. In particular, every finite group is a zero-dimensional Lie group. By a left \emph{$G$-action} on a manifold $X$, we mean a smooth map $G \times X \to X$ denoted by $(g,x) \mapsto g\cdot x$, where $(gg')\cdot x = g \cdot (g' \cdot x)$ and $e \cdot x = x$. Similarly, a right action is a map $(g,x)\mapsto x\cdot g$ with $x \cdot (gg') = (x\cdot g)\cdot g'$. In either case, we say that $G$ acts on $X$ on the left or right respectively. A \emph{$G$-space} is then a manifold $X$ along with a specified $G$-action.

\begin{definition}
We say that a vector bundle $E$ over $X$ is a \emph{$G$-vector bundle} if both $E$ and $X$ are $G$-spaces and the following two conditions are satisfied:
\begin{enumerate}
\item
The vector bundle projection map $\pi : E \to X$ commutes with the $G$-actions on $E$ and $X$. In other words, we have $g\cdot\pi(v) = \pi(g\cdot v)$ for all $v\in E$ and $g\in G$, and
\item
For each $g\in G$, the maps $E_{x}\to E_{g\cdot x}$ given by $v \mapsto g\cdot v$ are linear maps for all $x\in X$.
\end{enumerate}
\end{definition}

\begin{definition}
Given a $G$-space $X$, we define $K_{G}(X)$ by applying the Grothendieck construction to $\Vect_{G}(X)$, the semigroup of $G$-isomorphism classes of complex $G$-vector bundles over $X$. The study of $K_{G}(X)$ is called \emph{equivariant $K$-theory}.
\end{definition}

All of the relevant results about ordinary $K$-theory presented in \S\ref{sec:1} carry over directly to the equivariant case, making the appropriate changes in the notation. The most notable exceptions are that the splitting principle holds only when $G$ is abelian, and that the proof of the equivariant Thom isomorphism theorem is rather more difficult. For a complete discussion of equivariant $K$-theory, see \cite{Seg68}.

\begin{definition}
Given a continuous group homomorphism $\varphi : H \to G$, any $G$-space $X$ can be made into an $H$-space by taking the action $h\cdot x = \varphi(h) \cdot x.$ Applying this construction to $G$-vector bundles over $X$, we thus obtain a map $\varphi^{*}: K_{G}(X) \to K_{H}(X)$.
\end{definition}
 
Note that if $G$ is the trivial group, then $K_{G}(X)$ reduces to the ring $K(X)$ from ordinary $K$-theory. Also note that if $X$ is a point, then a complex $G$-vector bundle over $X$ is simply a complex $G$-module. In this
case $K_{G}(\pt)$ becomes the representation ring $R(G)$ of $G$, obtained by applying the Grothendieck construction to the semigroup of finite-dimensional complex representation spaces of $G$. In general, considering the homomorphism $R(G) \to K_{G}(X)$ induced by the map of $X$ onto a point, we see that $K_{G}(X)$ is a module over $R(G)$. In one extreme case, we have

\begin{proposition}
If $G$ acts trivially on $X$ (i.e., $g\cdot x = x\,\forall g$) then $K_{G}(X) \cong K(X) \otimes R(G).$
\end{proposition}

For a proof, see \cite[\S 2.2]{Seg68} or \cite[\S 1.6]{Ati67} for a discussion of the case of finite groups. At the other extreme, we obtain

\begin{proposition}
If $G$ acts freely on $X$ (i.e., $g\cdot x = x \leftrightarrow g = e$), then $K_{G}(X) \cong K(X/G),$ where $X/G$ is the orbit space obtained by identifying $x$ with $g\cdot x$ for all $g\in G$.
\end{proposition}

\begin{proof}
Let $E$ be a $G$-vector bundle over $X$. Then we see that $E/G$ is an ordinary vector bundle over $X/G$ since the $G$-actions are then trivial.\footnote{Note that it is also necessary to check that $E/G$ is locally trivial.}  We thus obtain a homomorphism $K_{G}(X) \to K(X/G)$. Going in the other direction, consider the map induced by the canonical projection $p : X \to X/G$ of $X$ onto its orbit space. Given an ordinary vector bundle $F$ over $X/G$, then induced bundle $p^{*}F$ has total space given by $\bigl\{ (x,v) \in X \times F \,\big|\, p(x) = \pi(v) \bigr\}$. By giving it the $G$-action $g \cdot (x,v) = (g \cdot x, v)$, we can make $p^{*}F$ into a $G$-vector bundle over $X$, thereby obtaining an inverse homomorphism $K(X/G) \to K_{G}(X)$. Using the canonical vector bundle isomorphisms $E \to p^{*}(E/G)$ and $(p^{*}F)/G \to F$, our result follows immediately.
\end{proof}

\subsection{Homogeneous Spaces}

One particularly interesting class of manifolds from the viewpoint of equivariant $K$-theory are the homogeneous spaces. Letting $H$ be a closed subgroup of $G$, we consider the quotient space $G/H$ of left cosets with the quotient topology. The resulting space is then a smooth manifold (see \cite[p.~21]{Ada69}). If $G$ acts transitively on a manifold $X$ (i.e., for each $x,y\in X$, we have $g\cdot x = y$ for some $g\in G$), then we say that $X$ is a \emph{homogeneous space}. Furthermore, if we let $H$ be the isotropy or stabilizer subgroup (i.e., the subgroup fixing a point), then we have an isomorphism $X \cong G/H$. Some simple examples of homogeneous spaces are:

\begin{itemize}
\item
Any Lie group: $G \cong G/ \{e\}$,
\item
The sphere: $S^{n} \cong SO(n+1)/SO(n)$,
\item
The odd sphere: $S^{2n+1}\cong U(n+1)/U(n)$ or $S^{2n+1}\cong SU(n+1)/SU(n)$ for $n> 0$,
\item
The real Grassmannian: $G(n,k)\cong O(n+k)/ O(n) \times O(k),$
\item
Real projective space: $\R P^{n} \cong O(n+1) / O(n) \times O(1),$
\item
The complex Grassmannian: $G_{\C}(n,k)\cong U(n+k)/ U(n) \times U(k),$
\item
Complex projective space: $\C P^{n} \cong U(n+1) / U(n) \times U(1),$
\item
The lens spaces: $L(n,k) = S^{2n+1}/\Z_{k} \cong U(n+1) / U(n) \times \Z_{k}.$
\end{itemize}

By the above proposition, noting that $H$ acts freely on $G$ (by the canonical action $h\cdot g = hg$), we obtain $K(G/H) \cong K_{H}(G)$. However, for our purposes, we are more interested in $G$-vector bundles over $G/H$. Given any finite-dimensional $H$-module $M$, we can construct the product $G \times_{H}M$ given by $G \times M$ modulo the action of $H$ as follows:
$$ G \times_{H} M = G \times M / (g,x) \sim h\cdot(g,x) = (gh^{-1},hx).$$
Letting $\pi : G \times_{H}M \to G/H$ be the projection map $(g,x) \mapsto gH$, we see that it clearly commutes with the canonical $G$-actions $g'\cdot gH = g'gH$ and $g'\cdot (g,x) = (g'g,x)$ on $G/H$ and $G\times_{H}M$, and so $G\times_{H}M$ is a vector bundle over $G/H$ with fiber $M$. In fact, we will now show that we can obtain any $G$-vector bundle of $G/H$ using this construction.

\begin{proposition}
If $H$ is a closed subgroup of a  Lie group $G$, then $K_{G}(G/H) \cong R(H)$.
\end{proposition}

\begin{proof}
Suppose we have a $G$-vector bundle $E$ over $G/H$. Noting that the $G$-action on $E$ restricts to an $H$-action on $E_{H}$, we see that the fiber of $E$ at the coset $H$ is a finite-dimensional $H$-module $M$. Then, the $G$-action $G\times E \to E$ on $E$ restricts to a $G$-homomorphism $G\times M \to E$ invariant under the action of $H$ on $G\times M$, and so it induces a $G$-map $\alpha : G \times_{H}M \to E$. We want to show that $\alpha$ is a $G$-isomorphism. By the definition of a $G$-vector bundle, we recall that the map $E_{H}\to E_{gH}$ given by $x \mapsto g \cdot x$ is an isomorphism for each $g\in G$ (considering that it has an inverse $x\mapsto g^{-1}\cdot x$). Restriction $\alpha$ to the fiber at any coset $gH$, we see that $\alpha|_{gH}: gH \times_{H} M \to E_{gH}$ is precisely that isomorphism, and so it follows that $\alpha$ is a global vector bundle isomorphism.
\end{proof}

\begin{example}
Considering the sphere $S^{2}\cong SO(3)/SO(2)$, we see that
$$K_{SO(3)} (S^{2}) \cong R\bigl(SO(2)\bigr) \cong \Z[z,z^{-1}].$$
In particular, we note that this is much larger ring that $K(S^{2}) \cong \Z^{2}$.
\end{example}

Essentially, ``translating'' by the $G$-action on the homogeneous space $G/H$, we saw that the global behavior of an equivariant vector bundle is determined by its behavior at a single point. By an entirely analogous argument, we can then prove a slight extension of this result.

\begin{proposition}
If $X$ is a locally compact\footnote{Here, we are using $K_{G}$-theory with compact supports, defined analogously to the ordinary case.} Hausdorff $H$-space, then $K_{G}(G\times_{H}X) \cong K_{H}(X)$.
\end{proposition}

\begin{example}
The tangent bundle of $G/H$ is a real vector bundle and at the coset $H$ it has fiber $\g/\h$, where $\g$ and $\h$ are the Lie algebras of $G$ and $H$.\footnote{The Lie algebra $\g \cong TG_{e}$ of a Lie group $G$ is the tangent space to $G$ at the identity, with bracket given by the bracket of left-invariant vector fields.} Noting that $T(G/H)$  is a $G$-vector bundle, and that $\g/\h$ has a natural $H$-action (the action $\Ad : H \to \Aut(\g/\h)$ induced by the map $A_{h} : gH \mapsto hgH$), it follows that $T(G/H) \cong G \times_{H}(\g/\h)$. Then, we obtain the isomorphism
$$K_{G}\bigl(T(G/H)\bigr) \cong K_{G}\bigl(G\times_{H}(\g/\h) \bigr)
\cong K_{H}(\g/\h).$$
\end{example}

\subsection{The $G$-Index Theorem}

If $E$ is a $G$-vector bundle over a compact $G$-space $X$, then we can make the vector space $\Gamma E$ of smooth sections of $E$ into an (infinite-dimensional) $G$-module as follows:

\begin{definition}
If $s : X \to E$ is a smooth section (not necessarily commuting with the $G$-actions on $X$ and $E$), then we give $\Gamma E$ the $G$-action $(g\cdot s)(x) = g\cdot s(g^{-1}\cdot x).$
\end{definition}

Given complex $G$-vector bundles $E$ and $F$ over $X$, the issue that naturally arises is to consider elliptic operators $P: \Gamma E \to \Gamma F$ that commute with the $G$-actions on $\Gamma E$ and $\Gamma F$. If this is the case, then we say that $P$ is \emph{$G$-invariant}. Noting that $TX$ is a $G$-space\footnote{Letting $L_{g} : X \to X$ be the  map $x\mapsto g \cdot x$, and letting $dL_{\g} : TX \to TX$ be its derivative, we define the canonical $G$-action on $TX$ by taking $g \cdot (x,\xi) = dL_{g}(x,\xi)$.} and that $\pi^{*}E$ and $\pi^{*}F$ are both $G$-vector bundles over $TX$, we see that the symbol $\sigma(P):\pi^{*}E \to \pi^{*}F$ is a $G$-map, and so it defines an element of $K_{G}(TX)$. As before, in the equivariant case we can define the topological index to be the composition
$$\tind_{G}^{X} : K_{G}(TX) \xrightarrow{i_{!}} K_{G}(TM) \xrightarrow{j_{!}^{-1}} K_{G}(\pt) \cong R(G),$$
where $i : X \to M$ is an embedding of $X$ into a finite-dimensional real $G$-module $M$, and $j : \pt \to M$ is the inclusion of the origin.

Since the operator $P$ is elliptic, we know that $\Ker P$ and $\Coker P$ are finite-dimensional. Furthermore, since $P$ is $G$-invariant, we see that $\Ker P$ and $\Coker P$ are invariant under the actions of $G$ on $\Gamma E$ and $\Gamma F$. It follows that the $G$-actions on $\Gamma E$ and $\Gamma F$ induce $G$-actions on $\Ker P$ and $\Coker P$, and so $\Ker P$ and $\Coker P$ are finite dimensional $G$-modules. We can thus define
\begin{definition}
The \emph{$G$-index} of a $G$-invariant elliptic operator $P$ is given by
$$\Index_{G} = [\Ker P] - [\Coker P] \in R(G).\footnote{We note that the $G$-index is related to the ordinary index by the equation $\Index P = \dim \Index_{G} P$.}$$
\end{definition}
Using the $G$-index, we can then define the analytical index to be the homomorphism
$$\aind_{G}^{X} : K_{G}(TX) \to K_{G}(\pt) \cong R(G)$$
given by $\aind_{G}^{X} = \Index_{G} P$, where $P$ is a $G$-invariant elliptic operator with $\sigma(P) = x.$ The proof of the Atiyah-Singer Index Theorem given in \S\ref{sec:3} naturally extends \emph{mutatis mutandis} (changing $\Z$ to $R(G)$ and inserting $G$'s everywhere) to the equivariant case. Hence, we obtain $\aind_{G} = \tind_{G}$, and the $G$-index theorem follows.

\begin{theorem}[$G$-Index Theorem]
If $P : \Gamma E \to \Gamma F$ is a $G$-invariant elliptic (pseudo-differential) operator on a compact manifold $X$, then
$$\Index_{G}P = \tind_{G}^{X}\sigma(P).$$
\end{theorem}

\subsection{The Induced Representation}

In this section, we review some of the basic results of the representation theory of Lie groups, stated in the notation of the representation ring. For a more complete discussion of representation theory, see \cite{Ada69}. For the rest of this section, we will take $G$ to be a compact Lie group. We say that a $G$-module is \emph{irreducible} if it has no $G$-invariant proper subspaces, and we say that two $G$-modules $M$ and $N$ are \emph{equivalent} if there exists a $G$-isomorphism $M\to N$. Letting $\{w_{\alpha}\}$ be a complete set of irreducible inequivalent representations of $G$, we see that the representation ring $R(G)$ is the free abelian group generated by $\{ w_{\alpha} \}$. Defining the \emph{intertwining number} of two $G$-modules by
$$\tau (M,N) = \dim_{\C}\Hom_{G}(M,N),$$
we can extend it to an inner product $\langle \cdot, \cdot \rangle_{G}$ on $R(G)$. Then by Schur's Lemma, if $M$ and $N$ are both irreducible $G$-modules, we have
$$\langle M, N\rangle_{G}  = \begin{cases}
1 & \text{ if $M$ and $N$ are equivalent,} \\
0 & \text{ otherwise,}
\end{cases}$$
and it follows that the set $\{ w_{\alpha} \}$ is orthonormal with respect to this inner product.

Next, we formally define the infinite extension $\hat{R}(G)$ of the representation ring to be the set of all possibly infinite linear combinations $\sum_{\alpha} a_{\alpha}w_{a}$ with $a_{\alpha}\in \Z$. Obviously, we have a canonical inclusion $R(G) \subset \hat{R}(G)$. The inner product on $R(G)$ then extends to a pairing $R(G) \otimes \hat{R}(G) \to \Z$ given by
$$\Bigl\langle x, \sum_{\alpha}a_{\alpha}w_{\alpha} \Bigr\rangle_{G}
= \sum_{\alpha} a_{\alpha} \langle x,w_{\alpha} \rangle_{G},$$
for each $x\in R(G)$ and $\sum_{\alpha}a_{\alpha}w_{\alpha}\in \hat{R}(G)$ (we note that this summation is actually a finite sum since $\langle x,w_{\alpha}\rangle$ vanishes for almost all $\alpha$).

\begin{definition}
Given a continuous group homomorphism $\varphi : H \to G$, we define the \emph{formal induced representation} $\varphi_{*} : R(H) \to \hat{R}(G)$ to be the map
$$\varphi_{*} : x \mapsto \sum_{\alpha}\langle \varphi^{*}w_{\alpha},x \rangle_{H} w_{\alpha},$$
where $\varphi^{*} :R(G) \to R(H)$ is the induced map from equivariant $K$-theory.
\end{definition}

Expressing this definition in terms of the extended inner product, we see that the induced representation satisfies the adjoint relation $\langle y, \varphi_{*}x\rangle_{G} = \langle \varphi^{*}y,x \rangle_{H}$ for all $x\in R(G)$ and $y\in R(G)$. If $H$ is a closed subgroup of $G$, and $i : H \to G$ is the inclusion map, then we can construct a more geometrically intuitive alternate definition for the induced representation. In particular, given a finite-dimensional $H$-module $M$, we recall that the space $\Gamma(G \times_{H}M)$ of smooth sections of the vector bundle $G \times_{H}M$ over $G/H$ has the structure of an infinite-dimensional $G$-module. Extending this map to the representation rings, we obtain

\begin{definition}
Given an inclusion $i : H \to G$ of a closed subgroup $H$, we define the \emph{induced representation} $i_{*} : R(H) \to \hat{R}(G)$ to be the extension to $R(H)$ of the map
$$M \mapsto \Gamma ( G \times_{H} M),$$
where $M$ is a finite-dimensional $H$-module.
\end{definition}

In order to show that these two definitions are equivalent, we will need to prove the following theorem of Frobenius:

\begin{theorem}[Frobenius Reciprocity]
If $W$ is a $G$-module and $M$ is an $H$-module, then
$$\Hom_{G}(W,i_{*}M) \cong \Hom_{H}(i^{*}W, M),$$
where $i_{*}M = \Gamma(G \times_{H}M)$ and $i^{*}W$ is the restriction of $W$ to an $H$-module.
\end{theorem}

\begin{proof}
Suppose that we have a $G$-invariant homomorphism $F : W \to \Gamma(G\times_{H}M)$. Then, composing it with the map $\Gamma(G\times_{H}M)\to M$ given by $s \mapsto s(H)$, we obtain the corresponding $H$-homomorphism $f:W\to M$, noting that for $h\in H$ and $w\in W$ we have
\begin{align*}
f(h\cdot w) &= F(h\cdot w)(H) = \bigr(h \cdot F(W)\bigr)(H) \\
&= h \cdot \bigl( F(w)(h^{-1}H) \bigr) = h \cdot \bigl( F(w)(H) \bigr)
= h \cdot f(w).
\end{align*}
Conversely, given an $H$-invariant homomorphism $f: W \to M$, we consider the corresponding map  $F : W \to \Gamma (G\times_{H}M)$ defined for $gH \in G/H$
and $w\in W$ by
$$F(w)(gH) = \bigl( g, f(g^{-1}\cdot w) \bigr) \in G \times_{H} M.$$
Since $f : W \to M$ is $H$-invariant, we see that this map is well defined, and we note that
\begin{align*}
F(g' \cdot w)(gH) &= \bigl( g, f(g^{-1}g' \cdot w) \bigr)
= g' \cdot \bigl( g'^{-1}g, f\bigl( (g'{-1} g)^{-1}\cdot w \bigr) \bigr) \\
&= g' \cdot F(w)(g'^{-1}\cdot gH) = \bigl( g' \cdot F(w)\bigr) (gH)
\end{align*}
for all $g\in G$. Hence, $F : W \to \Gamma(G\times_{H}M)$ is a $G$-homomorphism. Noting that the two maps $f \mapsto F$ and $F \mapsto f$ are inverses of each other, our result follows immediately.
\end{proof}

\begin{corollary}
If $i : H \to G$ is the inclusion of a closed subgroup $H$ in $G$, then we have
$$\bigl[ \Gamma (G \times_{H}M ) \bigr]
= \sum_{\alpha} \bigl\langle i^{*}w_{\alpha},[M] \bigr\rangle_{H}w_{\alpha}\in \hat{R}(G)$$
for any $H$-module $M$, and so the two definitions of the induced representation agree.
\end{corollary}

\begin{proof}
By the Frobenius reciprocity theorem, we see that
\begin{align*}
\bigl[ \Gamma(G\times_{H}M) \bigr]
&= \sum_{\alpha}\dim\Hom_{G}\bigl( w_{\alpha},\Gamma(G\times_{H}M) \bigr) \cdot w_{\alpha} \\
&= \sum_{\alpha}\dim\Hom_{H}(i^{*}w_{\alpha}, M) \cdot w_{\alpha}
= \sum_{\alpha} \bigl\langle i^{*}w_{\alpha},[M] \bigr\rangle_{H}w_{\alpha},
\end{align*}
as we wanted to show.
\end{proof}

\subsection{Homogeneous Differential Operators}

As a special case of the $G$-index theorem, we now consider the $G$-index of a $G$-invariant pseudo-differential\footnote{All of the results in this section apply equally well to pseudo-differential operators, but for notational convenience, we will omit the ``pseudo-'' prefix for the remainder of the section.} operator $P : \Gamma E \to \Gamma F$ on a homogeneous space $G/H$. If this is the case, then we say that $P$ is a \emph{homogeneous differential operator}. In light of the Peter-Weyl theorem, the analysis involved in this case becomes the representation theory of the induced map.  Also, as a result of the excessive structure imposed on the class of homogeneous operators, we will see shortly that the index of $P$ depends only on the vector bundles $E$ and $F$ and not directly on the operator $P$ itself. Following \cite{Bott65} we define the homogeneous symbol by

\begin{definition}
Given a homogeneous operator $P : \Gamma E \to \Gamma F$ over $G/H$, we recall that we can write $E$ and $F$ in the form $E \cong G \times_{H}M$ and $F\cong G\times_{H}N$, where  $M$ and $N$ are finite-dimensional $H$-modules. We thus define the \emph{homogeneous symbol} of $P$ to be
$$\hat{\sigma}(P) = [M] - [N] \in R(H).$$
\end{definition}

We now compare this homogeneous symbol to the ordinary symbol $\sigma(P)\in K_{G}\bigl(T(G/H)\bigr),$ represented by the homogeneous complex
$$\xymatrix{
0 \ar[r] & \pi^{*}E \ar[r]^{\sigma} & \pi^{*}F \ar[r] & 0
}$$
over $T(G/H)$, where $\sigma : \pi^{*}E \to \pi^{*}F$ is a $G$-homomorphism. Letting $s : \pt \to T(G/H)$ denote the inclusion of the identity coset $H\in G/H$ into $T(G/H)$ via the zero-section, we see that the induced map
$s^{*}: K_{G}\bigl(T(G/H)\bigr) \to K_{G}(\pt) = R(G)$ maps $\sigma(P)$ to the complex
$$\xymatrix{
0 \ar[r] & M \ar[r]^{\sigma(H,0)} & N \ar[r] & 0
}$$
over a point.\footnote{Note that the map $\sigma(H,0): M \to N$ is a homomorphism, and is generally \emph{not} an $H$-homomorphism.} It then follows that
$$ s^{*}\sigma(P) = [M] - [N] = \hat{\sigma}(P) \in R(H).$$
We can now compute the $G$-index of a homogeneous differential operator.

\begin{theorem}[Bott]
If $P$ is a homogeneous differential operator over $G/H$, then
$$\Index_{G}P = i_{*}\hat{\sigma}(P) = i_{*} \circ s^{*}\sigma(P) \in \hat{R}(G),$$
where $i_{*}: R(H) \to \hat{R}(G)$ is the formal induction homomorphism, and in particular, $i_{*}\hat{\sigma}(P)$ is actually a finite element contained in $R(G) \subset \hat{R}(G)$.
\end{theorem}

By our alternate definition of the induced representation $i_{*} : R(H) \to \hat{R}(G)$, we have
$$i_{*}\bigl( \hat{\sigma}(P) \bigr)
= i_{*}\bigl( [M] - [N] \bigr)
= \bigl[ \Gamma (G\times_{H}M ) \bigr] - \bigl[ \Gamma(G\times_{H}N) \bigr]
= [\Gamma E] - [\Gamma F ],$$
and so we can restate Bott's theorem even more concisely.

\begin{theorem}
If $P : \Gamma E \to \Gamma F$ is a homogeneous differential operator, then
$$\Index_{G}P = [\Gamma E] - [\Gamma F] \in \hat{R}(G),$$
where $[\Gamma E] -[\Gamma F]$ is in fact a finite element contained in $R(G) \subset \hat{R}(G)$.
\end{theorem}

In order to prove this theorem, we need to find a decomposition of the infinite-dimensional $G$-module $\Gamma E$ into finite-dimensional $G$-invariant subspaces. Restricting ourselves to such subspaces, the theorem can be much more readily verified, and we can then lift the result to the global space of smooth sections. The natural candidates for these subspaces are defined by

\begin{definition}
For each $w_{\alpha}$, we define $\Gamma_{\alpha}E$ to be the image of the homomrorphism
$$i_{\alpha} : w_{\alpha}\otimes \Hom_{G}(w_{\alpha},\Gamma E) \to \Gamma E$$
given by $i_{\alpha}: (w \otimes \varphi) \mapsto \varphi(w)$. We then see that $\Gamma_{\alpha}(E)$ is a finite dimensional $G$-invariant subspace of $\Gamma E$, with image in the representation ring $[\Gamma_{\alpha}E] = \langle w_{\alpha},[\Gamma E] \rangle_{G} w_{\alpha}\in R(G)$.
\end{definition}

Noting that the class of $\Gamma E$ in the infinite extensions $\hat{R}(G)$ of the representation ring is given by $[\Gamma E] = \sum_{\alpha}\langle w_{\alpha},[\Gamma E]\rangle_{G}w_{\alpha}$, this expression suggests that we can (at least symbolically) obtain an infinite decomposition of the form $\Gamma E = \oplus_{\alpha}\Gamma_{\alpha} E.$ Such a decomposition is made explicit by the following famous result:

\begin{lemma}[Peter-Weyl Theorem]
If we let $L^{2}(E)$ be the Hilbert space obtained by taking the completion of the space $\Gamma E$ with respect to a $G$-invariant inner product, then the subspaces $\Gamma_{\alpha}E$ form a complete set of orthogonal subspaces of $L^{2}(E)$.
\end{lemma}

We now have enough background to give a fairly quick proof of Bott's theorem.

\begin{proof}[Proof of Bott's Theorem]
Let $P : \Gamma E \to \Gamma F$ be an elliptic homogeneous differential operator. Since $P$ is a $G$-map, we see that $P$ induces linear maps
$P_{\alpha}:\Gamma_{\alpha}E \to \Gamma_{\alpha} F$ on each of the finite-dimensional $G$-invariant subspace $\Gamma_{\alpha}$. Since $P$ is elliptic, we know that $\Ker P$ and $\Coker P$ are finite-dimensional $G$-modules, and the Peter-Weyl theorem tells us that
$$\Ker P = \bigoplus_{\alpha}\Ker P_{\alpha}, \qquad
  \Coker P = \bigoplus_{\alpha}\Coker P_{\alpha},$$
where $\Ker  P_{\alpha}$ and $\Coker P_{\alpha}$ are trivial for almost all $\alpha$. We thus obtain
$$\Index_{G}P = [\Ker P] - [\Coker P] = \sum_{\alpha}[\Ker P_{\alpha}] - [\Coker P_{\alpha}].$$
Then, noting that we have the exact sequence
$$\xymatrix{
	0 \ar[r]
	& \Ker P_{\alpha} \ar[r] 
	& \Gamma_{\alpha}E \ar[r]^{P_{\alpha}}
	& \Gamma_{\alpha}F \ar[r] 
	& \Coker P_{\alpha} \ar[r]
	& 0
}$$
of finite-dimensional $G$-modules, we can apply the $K_{G}$-theoretic construction of $R(G)$ in terms of complexes. Since this complex has empty support, its Euler characteristic (i.e., the alternating sum of the $G$-modules)
in $R(G)$ must vanish. This then gives us
\begin{align*}
	[\Ker P_{\alpha}] - [\Coker P_{\alpha}]
	&= [\Gamma_{\alpha}E] - [\Gamma_{\alpha}F] \\
	&= \bigl\langle w_{\alpha}, [\Gamma E] \bigr\rangle_{G}w_{\alpha}
	= \bigl\langle w_{\alpha}, [\Gamma F] \bigr\rangle_{G} w_{\alpha.}
\end{align*}
Summing over $\alpha$ we obtain
$$\Index_{G}P
= \sum_{\alpha}\bigl\langle w_{\alpha},[\Gamma E] \bigr\rangle_{G}w_{\alpha}
- \sum_{\alpha}\bigl\langle w_{\alpha},[\Gamma F] \bigr\rangle_{G}w_{\alpha}
= [\Gamma E] - [\Gamma F],$$
our desired result\end{proof}

\begin{example}
Consider the case of the circle $S^{1}$, taking $G = SO(2)$ and letting $H$ be the trivial group. Then we see that any elliptic homogeneous pseudo-differential operator on the circle must be of the form $P : \Gamma(S^{1}\times \C^{n}) \to \Gamma (S^{1}\times\C^{n})$. In this case, the homogeneous symbol $\hat{\sigma}(P) = [\mathbf{n}] - [\mathbf{n}] = 0$ vanishes, and so we have $\Index_{G}P = 0$. Furthermore, since $P$ is $G$-invariant, we see that $P$ must be given by $P = p(-i\,d/d\theta)$, where $p(\xi)$ is an $n\times n$ matrix-valued function with constant coefficients.

It is easy to see that the inequivalent irreducible representations of $S^{1} = SO(2)$ are the one-dimensional representations given by $w_{m}: \theta \to e^{im\theta}$ for $m\in \Z$. We then see that the finite-dimensional subspaces of $\Gamma(S^{1}\times\C^{n})$ are given by $\Gamma_{m}(S^{1}\times \C^{n}) = nw_{m}$, and so the Peter-Weyl theorem tells us that any smooth function $f : S^{1}\to \C^{n}$ can be written as
$$ f = \sum_{m}a_{m}w_{m} : \theta \mapsto \sum_{m} a_{m}e^{im\theta},$$
where $a_{m}\in \C^{n}$. This is just the Fourier series expansion. The restrictions
$$P_{m} : \Gamma_{m}(S^{1}\times \C^{n}) \to \Gamma_{m}(S^{1}\times\C^{n})$$
are then given by
$$P_{m}(\theta \mapsto a_{m}e^{im\theta}) = \theta \mapsto p(m) a_{m}e^{im\theta},$$
and we see that $\Ker P_{m} = \Ker P(m)$ and $\Coker P_{m} = \C^{n}/ \im P(m) \cong \Ker P(m)$. It follows that $[\Ker P_{m}] - [\Coker P_{m}] =0$, which extends to the global result that $\Index_{G}P = 0$.
\end{example}

\appendix

\section{A Proof of Bott Periodicity}
\label{sec:appendix}

\subsection{Preliminaries}

Throughout the main body of this paper, we have been using the machinery of cohomology, $K$-theory, and representation theory to solve the index problem from analysis. In this appendix, the tables are turned. We will now present an application of index theory to topology, in particular using the index of a family of elliptic operators to give a simple proof of the Bott periodicity theorem. We recall from \S\ref{sec:1} that the Bott periodicity theorem provides the basic foundation of $K$-theory, giving it the structure of a periodic cohomology theory. Hopefully, this proof will impress the reader with the fundamental connection between algebraic topology and elliptic operators (if the reader has not already been so impressed by the remainder of this paper). For more discussion of this connection, see \cite{Ati67a}. Our proof will be along the lines of \cite{Ati69}, presenting a slightly simpler case of the more general argument given in \cite{Ati68}. Later in this appendix, we will sketch how this general argument can be extended to prove the equivariant version of the Thom isomorphism theorem.

Before proceeding, we recall our statement of the Bott periodicity theorem from \S\ref{sec:1}. Let $b$ be a generator of $\tilde{K}(S^{2}) = K(\C)$. Noting that the exterior tensor product $\tensorhat : K(X) \otimes K(Y) \to K(X\times Y)$ restricts to the product $\tensorhat : \tilde{K}(X) \otimes \tilde{K}(Y) \to \tilde{K}(X \wedge Y),$\footnote{Using the decomposition $K(X) \cong \tilde{K}(X)\oplus K(\pt)$, the exterior tensor product then splits as
$$\otimes : \bigl( \tilde{K}(X) \otimes \tilde{K}(Y) \bigr) \oplus
\bigl( \tilde{K}(X) \oplus \tilde{K}(Y) \bigr)  \oplus K(\pt)
\to \tilde{K}(X\wedge Y) \oplus \tilde{K}(X \vee Y) \oplus K(\pt).$$
Since $\tilde{K}(X)\oplus\tilde{K}(Y) \to \tilde{K}(X\vee Y)$, we see that the exterior tensor product restricts to the product map
$$\tensorhat : \tilde{K}(X) \otimes \tilde{K}(Y) \to \tilde{K}( S\wedge Y),$$
where we recall that $X\vee Y = (X \times \pt) \cup (\pt \times Y)$ is the one-point union, and $X\wedge Y = ( X\times Y ) / (X\vee Y)$ is the smash product.}
we define the homomorphism
$$\beta : K(X) \to K^{-2}(X) = \tilde{K}(S^{2}\wedge X^{+})$$
by $\beta(X) = b \tensorhat x$. The Bott periodicity theorem states that this map is an isomorphism. Using $K$-theory with compact supports, we note that
$$K(\R^{n}\times X) = \tilde{K}\bigl( (\R^{n}\times X)^{+}\bigr)
= \tilde{K}(S^{n}\wedge X^{+}) = K^{-n}(X),$$
and so we can view $\beta$ as a homomorphism
$$\beta : K(X) \xrightarrow{b \tensorhat \cdot} K(\R^{2}\times X),$$
where $b$ is viewed as an element of $K(\R^{2}).$ In order to prove the Bott periodicity theorem, we will use the index of a family of elliptic operators to construct an inverse homomorphism $\alpha : K(\R^{2}\times X) \to K(X)$.

\subsection{The Wiener-Hopf Operator}

We begin by considering a specific pseudo-differential operator on the circle $S^{1} = \R/2\pi\Z$. Recall from \S\ref{sec:3} that we introduced the pseudo-differential operator $P$ of order zero on complex-valued functions on $S^{1}$ given in terms of Fourier series by
$$P(e^{ik\theta}) = \begin{cases}
e^{ik\theta} & \text{ for $k\geq 0$,} \\
0 & \text{ for $k < 0$.}
\end{cases}$$
This operator can be viewed as a projection from the space of complex-valued functions on $S^{1}$ to the subspace consisting of functions with only positive Fourier coefficients. Given any continuous function $f : S^{1}\to \C^{*}$, where $\C^{*}$ is the non-zero complex numbers, we can construct a corresponding pseudo-differential operator $P_{f} = P\,f\, P + (1-P)$. Writing $f(\theta) = \sum_{n}a_{n}e^{in\theta}$, we see that in terms of the Fourier series expansion, we obtain
$$P_{f}(e^{ik\theta}) = \begin{cases}
\sum_{n>-k}a_{n}e^{i(n+k)\theta} & \text{ for $k\geq 0$,} \\
e^{ik\theta} & \text{ for $k < 0$,}
\end{cases}$$
and we note that $P_{f}$ is elliptic since its $0$-th order symbol is given by
$$\sigma_{P_{f}}(\theta,\xi) = \begin{cases}
f(\theta) & \text{ for $\xi > 0$,} \\
1 & \text{ for $\xi < 0$.}
\end{cases}$$
In particular, defining $P_{n} = P_{f}$ where $f(\theta) = e^{in\theta}$,
we see that for $n\geq 0$ we have
$$P_{n}(e^{ik\theta}) = \begin{cases}
e^{i(n+k)\theta} & \text{ for $k \geq 0$,} \\
e^{ik\theta} & \text{ for $k < 0$,}
\end{cases}$$
and so $\Ker P_{n}$ is trivial while $\Coker P_{n}$ is the $n$-dimensional space of functions of the form $g(\theta) = \sum_{0\leq m < n}b_{m}e^{imx}$. Hence, $\Index P_{n} = -n$. Similarly, for $n<0$ we have
$$P_{n}(e^{ik\theta}) = \begin{cases}
e^{i(n+k)\theta} & \text{ for $k \geq -n$,} \\
0 & \text{ for $0 \leq k < -n$,} \\
e^{ik\theta} & \text{ for $k < 0$.}
\end{cases}$$
In this case $\Coker P_{n}$ is trivial while $\Ker P_{n} \cong \Coker P_{-n}$, and so we again obtain the result $\Index P_{n} = -n$. Now, we recall that any continuous map $f : S^{1}\to \C^{*}$ is homotopic to a map of the form $f(\theta) = e^{in\theta}$, where we define the degree of the map to be $\deg f = n$. Since the index is a locally constant function, we then see that $\Index P_{f} = -\deg f$.

We now generalize this construction. First, we can without loss of generality extend our discussion to complex vector-valued functions on $S^{1}$. Given a continuous function $f : S^{1}\to GL(N,\C)$, we can then construct the corresponding elliptic operator $P_{f}$, and we note that $\Index P_{f}$ depends only on the homotopy class of $f$, which in this case is the class $[f] \in \pi_{1}(GL(N,\C))$ in the fundamental group of $GL(N,\C)$. Taking this one step further, let $V$ be a complex vector bundle over a compact space $X$, and suppose that we have a continuously varying collection of continuous maps $f_{x} : S^{1}\to \Aut V_{x}$ for each $x\in X$. We then obtain a family
$P_{f}$ of pseudo-differential operators on $S^{1}$ parametrized by $X$, where each of the operators $(P_{f})_{x} = P_{f_{x}}$ is elliptic. Such a collection is called an \emph{elliptic family}. If $\dim\Ker (P_{f})_{x}$ and $\dim\Coker(P_{f})_{x}$ are constant for all $x\in X$, then we can define $\Ker P_{f}$ and $\Coker P_{f}$ to be the corresponding vector bundles over $X$, and so we can define
$$\Index P_{f} = [\Ker P_{f}] - [\Coker P_{f}] \in K(X).$$
For the general construction of the index of a family of Fredholm operators in the case where $\dim\Ker (P_{f})_{x}$ and $\dim\Coker (P_{f})_{x}$ are not constant, see \cite[\S 2]{Ati69}. We note that this index depends only on the homotopy class $[f]$ of the map $f : (x,\theta) \mapsto \Aut V_{x}.$

\subsection{The Clutching Construction}

Given a vector bundle $V$ over $X$ and a continuous map of the form $f : (x,\theta) \mapsto \Aut V_{x}$, we can construct a vector bundle $E$ over $S^{2}\times X$ as follows. Letting $B^{+}$ and $B^{-}$ be the closed upper and lower hemispheres of the sphere $S^{2}$, we see that their union $B^{+}\cup B^{-} = S^{2}$ is the whole sphere, while their intersection $B^{+}\cap B^{-} = S^{1}$ is the equator. We then obtain $E$ by taking the union relative to the map $f$
$$E = (B^{-}\times V) \cup_{f} (B^{+}\times V),$$
where along the intersection $(B^{-}\times V)\cap (B^{+}\times V) = S^{1}\times V$ we have the identification
$$(\theta, v) \in B^{-}\times V_{x} \quad \sim \quad
\bigl( \theta, f(\theta,x) v \bigr) \in B^{+}\times V_{x}.$$
This construction is known as the \emph{clutching construction} and the map $f : (x,\theta) \mapsto \Aut V_{x}$ is called a \emph{clutching function}. For a complete discussion of clutching functions, see \cite[\S 1.4]{Ati67}.

We note that every vector bundle $E$ over $S^{2}\times X$ can be obtained by this construction\footnote{Letting $i_{+},i_{-} : X \to S^{2}\times X$ be the inclusions at the two poles of $S^{2}$, we note that since $B^{+}$ and $B^{-}$ are contractible, we have $E|_{B^{+}\times X}\cong B^{+}\times i_{+}^{*}E$ and $E|_{B^{-}\times X}\cong B^{+}\times i_{-}^{*}E$. Putting $V \cong i^{*}_{-}E \cong i^{*}_{+}E$, we then take $f : (x,\theta) \mapsto \Aut V_{x}$ to be the identification of $E|_{B^{+}\times X}$ and $E|_{B^{-}\times X}$ on $S^{1}\times X$.} from a vector bundle $V$ over $X$ and an appropriate clutching function $f : (x,\theta)\mapsto \Aut V_{x}$, and in addition, $E$ is uniquely determined (up to isomorphism) by the homotopy class $[f]$ of $f$ (see \cite[\S 1.4.6]{Ati67}). We thus obtain a map $\Vect(S^{2}\times X) \to K(X)$ given by
$$ E \mapsto (V,[f]) \mapsto P_{f}\mapsto \Index P_{f}\in K(X),$$
which then extends to a group homomorphism $\alpha : K(S^{2}\times X) \to K(X)$.
Noting that\footnote{Considering the sequence of spaces $X^{+}\vee Y^{+}\to X^{+}\times Y^{+}\to X^{+}\wedge Y^{+} = (X\times Y)^{+}$, we obtain the exact sequence $0 \to K(X\times Y) \to \tilde{K}(X^{+}\times Y^{+}) \to K(X) \oplus K(Y) \to 0$. In particular, taking $X = \R^{n}$ and $Y = X$, we have the exact sequence $0 \to K(\R^{n}\times X)\to \tilde{K}(S^{n}\times X^{+}) \to \tilde{K}(S^{n}) \oplus K(X) \to 0$. Noting that $\tilde{K}(S^{n}\times X^{+}) = \tilde{K}(S^{n}) \oplus K(S^{n}\times X)$, we see that $K(S^{n}\times X) \cong K(\R^{n}\times X) \oplus K(X)$.} $K(S^{2}\times X) \cong K(\R^{2}\times X) \oplus K(X)$, we can take the restriction of $\alpha$ to $K(\R^{2}\times X)$ to obtain
$$\alpha : K(\R^{2}\times X) \to K(X).$$

Taking $X$ to be a point, a vector bundle over $X$ is simply a vector space, and for one-dimensional bundles, a clutching function $f : (x,\theta) \mapsto \Aut V_{x}$ is simply a continuous map $f : S^{1}\to \C^{*}$. Then for each $n\in\Z$ we define the bundle $E_{n}$ over $S^{2}$ by applying the clutching construction with $V = \C$ and $f(\theta) = e^{in\theta}$. We note that it follows immediately from our construction that $\alpha[E_{n}] = -n$. Recalling that every such map $f$ is homotopic to one of the form $f(\theta) = e^{in\theta}$, we see that this construction yields all of the one-dimensional bundles over $S^{2}$. In particular, this tells us that $\tilde{K}(S^{2})\cong \Z$. We will now look at several examples of vector bundles over $S^{2}$ from this point of view.

\begin{example}
Viewing $S^{2}\cong \C\cup\infty$ as the Riemann sphere, we can consider its complex tangent bundle $T_{\C}S^{2}$. We can then identify $B^{-}$ with the unit disc $\bigl\{ z \, \big| \, |z| \leq 1 \bigr\}$, and we identify $B^{+}$ with $\bigl\{ w \, \big| \, |w| \leq 1 \bigr\}$, where $w = z^{-1}$ (since we are dealing with the complex structure, we want our coordinate transitions to be complex analytic). Now, a tangent vector at a point $z\in B^{-}$ can be expressed in the form $\xi = d/dt\,f(t)|_{t=0}$, where $f(0) = z$. Translating to $w$ coordinates, the corresponding tangent vector is given by
$$\xi_{w} = \frac{d}{dt}\,f(t)^{-1} \Bigr|_{t=0}
= -f(t)^{-2}\,\frac{d}{dt}\,f(t)\Bigr|_{t=0} = -z^{-2}\xi_{z}.$$
We can thus write $$T_{\C}S^{2} \cong \bigl( \,( B^{-}\times \C) \cup (B^{+}\times \C) \, \bigr) \,/\, \sim,$$ where we have the identification on the unit circle $\bigl\{ z \, \big| \, |z| = 1 \bigr\} = B^{-}\cap B^{+}$ given by
$$(z,\xi)\in B^{-}\times \C \; \sim \; (z,-z^{-2}\xi) \in B^{+}\times \C.$$
By rotating halfway around the unit circle, we see that $f(z) = -z^{-2}$ is homotopic to the map $f(x) = z^{-2}$. Since this map corresponds to $f(\theta) = e^{-2i\theta}$, we see that $T_{\C}S^{2}\cong E_{-2}$.
\end{example}

\begin{example}
Considering the map $f(\theta) = e^{in\theta}$, we see that the symbol $\sigma(P_{n})\in K(TS^{1})$ of the corresponding elliptic pseudo-differential operator $P_{n}$ is given by the complex
$$\xymatrix{
0 \ar[r] & TS^{1}\times \C \ar[r]^{\sigma_{n}} & TS^{1}\times \C \ar[r] & 0}$$
over $TS^{1}$, where the homomorphism is given by
$$\sigma_{n}(\theta,\xi) = \begin{cases}
e^{in\theta} & \text{ for $\xi > 0$,} \\
1 & \text{ for $\xi < 0$.}
\end{cases}$$
Letting $p^{+}\in B^{+}$ and $p^{-}\in B^{-}$ be the ``north'' and ``south'' poles of the sphere $S^{2}$, we can identify $B^{-}-\{p^{-}\}$ with
$\bigl\{ ( \theta,\xi)\in TS^{1} \,\big|\, \xi \leq 0 \bigr\}$ and $B^{+}-\{p^{+}\}$ with
$\bigl\{ (\theta,\xi)\in TS^{1}\,\big|\,\xi \geq 0 \bigr\}$. We note that $\sigma_{n}(\theta,\xi)$ is homotopic to the homomorphism $\sigma'_{n}(\theta,\xi)$ given by
$$\sigma'(\theta,\xi) = \begin{cases}
f(\theta) / (1 + \xi) & \text{ for $\xi > 0$,} \\
1 & \text{ for $\xi < 0$,}
\end{cases}$$
which we can extend to a continuous map $\sigma'_{n} : E_{0}\to E_{n}$ by defining
$$\sigma'_{n}(p^{-}) = \lim_{\xi\to-\infty}\sigma'_{n}(\theta,\xi) = 1,
\qquad
\sigma'_{n}(p^{+}) = \lim_{\xi\to+\infty}\sigma'_{n}(\theta,\xi) = 0.$$
We note that this map is continuous on the circle $S^{1}\times\{0\}$ since by our clutching construction we are identifying the point $\bigl( (\theta,0),z\bigr) \in B^{-}\times \C$ with $\bigl( (\theta,0),e^{in\theta}z \bigr) \in B^{+}\times \C$. Hence the symbol complex extends to a complex
$$\xymatrix{
0 \ar[r] & E_{0} \ar[r]^{\sigma'_{n}} & E_{n} \ar[r] & 0}$$
over $S^{2}$, and we have $\sigma(P_{n}) = [E_{0}] - [E_{n}] \in \tilde{K}(S^{2}).$ It follows that $\alpha \bigl(\sigma(P_{f})\bigr) =n$.
\end{example}

\begin{proposition}
$\aind_{S^{1}}(x) = -\alpha_{\pt}(x)$ for any $x\in K(TS^{1})\cong K(\R^{2}) \cong \tilde{K}(S^{2}).$
\end{proposition}

\begin{proof}
Noting that $K(TS^{1}) \cong \Z$, we see that every element $x\in K(TS^{1})$ occurs as the symbol $x\in \sigma(P_{n})$ of one of the elliptic pseudo-differential operators $P_{n}$. We then have
\begin{align*}
\aind_{S^{1}}(x) &= \Index P_{n} = -n \\
&= -\alpha_{\pt}\bigl( \sigma(P_{n}) \bigr) = -\alpha_{\pt}(x)
\end{align*}
by the above discussion.
\end{proof}

\begin{example}
We recall from \S\ref{sec:1} that the Bott class $\lambda_{1}\in K(\C)\cong \tilde{K}(S^{2})$ is represented by the exterior complex $\Lambda(\C)$ over $\C$ given by
$$\xymatrix{
0 \ar[r] & \pi^{*}\Lambda^{0}(\C) \ar[r]^{\alpha} & \pi^{*}\Lambda^{1}(\C) \ar[r] & 0 },$$
where $\pi^{*}\Lambda^{0}(\C) \cong \pi^{*}\Lambda^{1}(\C) \cong \C \times \C$ and $\alpha : (z,v) \mapsto (z,zv)$. Writing $D^{-} = \bigl\{ z\,\big|\,|z| \leq 1\bigr\}$ and $D^{+} = \bigl\{ z\,\big|\,|z| \geq 1\bigr\}$, we can identify $D^{-}$ with $B^{-}$ and $D^{+}\cup \infty$ with $B^{+}$ via the inclusion $\C \subset \C^{+}\cong S^{2}$.
Then, since $E_{0}$ is the trivial bundle over $S^{2}$, we see that $\pi^{*}\Lambda^{0}(\C) = E_{0}|_{\C}$, and we can construct an isomorphism $\beta : \pi^{*}\Lambda^{1}(\C) \to E_{-1}|_{\C}$ by taking
$$\beta : (z,v) \mapsto \begin{cases}
(z,v)\in B^{-}\times \C & \text{ for $z\in D^{-}$,} \\
(z,z^{-1}v) \in B^{+}\times \C & \text{ for $z\in D^{+}$.}
\end{cases}$$
Identifying the circle $S^{1} = \R / 2\pi\Z$ with the unit circle in $\C$, we note that the map $f(\theta) = e^{-i\theta}$ is equivalent to the map $f(z) = z^{-1}$, and so $\beta$ is continuous by our clutching construction. The exterior complex $\Lambda(\C)$ can thus be rewritten as the complex
$$\xymatrix{
0 \ar[r] & E_{0}|_{\C} \ar[r]^-{\gamma} & E_{-1}|_{\C} \ar[r] & 0}$$
over $\C$, where $\gamma = \beta \circ \alpha$ is given by
$$\gamma : (z,v) \mapsto \begin{cases}
(z,zv) \in B^{-}\times \C & \text{ for $z\in D^{-}$, } \\
(z,v) \in B^{+}\times \C & \text{ for $z\in D^{+}$. }
\end{cases}$$
Finally, taking $\gamma(\infty,v) = (\infty,v)$, we can extend this to the complex
$$\xymatrix{
0 \ar[r] & E_{0} \ar[r]^-{\gamma} & E_{-1} \ar[r] & 0}$$
over $S^{2}$, and since $S^{2}$ is compact it follows that $\lambda_{1} = [E_{0}] - [E_{-1}]\in \tilde{K}(S^{2}) = K(\R^{2}).$ Computing that $\alpha(\lambda_{1}) = \alpha[E_{0}] - \alpha[E_{-1}] = 0-1 = -1$, the above proposition tells us that $\aind_{S^{1}}(\lambda_{1}) = -\alpha(\lambda_{1}) = 1$, which we needed in \S\ref{sec:3} to complete our proof of the normalization axiom.
\end{example}

\subsection{The Bott Periodicity Theorem}

For our generator of $K(\R^{2})$, we take $b = -\lambda_{1} = [E_{-1}] -[E_{0}],$ and so we obtain $\alpha(b) = 1$. We also note that from our construction we have the commutative diagram
$$\xymatrix{
	K(\R^{2}\times X)\otimes K(Y) \ar[r]^-{\tensorhat} \ar[d]^{\alpha_{x}\otimes 1}
	& K(\R^{2}\times X\times Y ) \ar[d]^{\alpha_{X \times Y}} \\
	K(X) \otimes K(Y) \ar[r]^-{\tensorhat} & K(X \times Y)
}$$
which tells us that for any $x\in K(\R^{2}\times X)$ and $y\in K(Y)$ we have the product formula
$$\alpha(X) \cdot y = \alpha (x\cdot y).$$
We also note that our map $\alpha : K(\R^{2}\times X) \to K(X)$ can be extended to locally compact spaces $X$ by taking the one point compactification $X^{+}$, and the product formula still holds. Now, the Bott periodicity theorem follows directly from a few simple algebraic properties of $K$-theory.

\begin{theorem}[Bott Periodicity]
Given any locally compact space $X$, then the two homomorphisms
\begin{align*}
\beta &: K(X) \to K(\R^{2}\times X) \\
\alpha &: K(\R^{2}\times X) \to K(X)
\end{align*}
are inverses of each other.
\end{theorem}

\begin{proof}
First, by our product formula, we note that for any $x\in K(X)$ we have
$$\alpha\beta(X) = \alpha(b\cdot x) = \alpha(b) \cdot x  = 1 \cdot x = x,$$
and so $\alpha$ is a left inverse of $\beta$. We must now show that $\alpha$ is also a right inverse of $\beta$ by verifying that for each $u \in K(\R^{2}\times X)$ we have $b\cdot \alpha(u) = \beta\alpha(u) = u$. Noting that our product formula gives us $\alpha (u \cdot b ) = \alpha(u) \cdot b$,
where $u\cdot \lambda_{1}\in K(\R^{2}\times X \times \R^{2})$, we would like to commute these two products. To do this, we consider the map
$$\tau : \R^{2}\times X \times \R^{2} \to \R^{2}\times X \times \R^{2}$$
that interchanges the two copies of $\R^{2}$. Taking $X$ to be a point, the corresponding transformation $\tau : \R^{4}\to \R^{4}$ is given by the matrix
$$ \begin{pmatrix} 0 & I_{2} \\ I_{2} & 0 \end{pmatrix}, \text{ where }
I_{2} = \begin{pmatrix} 1 & 0 \\ 0 & 1 \end{pmatrix},$$
which has determinant $+1$, and so it is in the same connected component of $GL(4,\R)$ as the identity. For general $X$, we then see that we can construct a homotopy between $\tau$ and the identity, and so the induced map of $K$-theory
$$\tau^{*}: K(\R^{2}\times X \times \R^{2}) \to K(\R^{2}\times X \times \R^{2})$$
is the identity map.  Letting $\tilde{u}\in K(X\times \R^{2})$ be the element corresponding to $u$ under the canonical isomorphism $K(\R^{2} \times X) \to K(X\times \R^{2})$, we obtain
$$u \cdot b = \tau^{*}(u \cdot b) = b \cdot \tilde{u},$$
and applying $\alpha$ gives us
$$ \alpha(u) \cdot b = \alpha(u\cdot b) = \alpha(b\cdot \tilde{u}) = \tilde{u} \in K(X \times \R^{n}).$$
Transforming back via the isomorphism $K(X\times \R^{2}) \to K(\R^{2}\times X)$, we  conclude that
$$u = b \cdot \alpha (u) = \beta\alpha(u) \in K(\R^{2}\times X),$$
which completes our proof.
\end{proof}

\subsection{The Thom Isomorphism Theorem}

This proof of Bott periodicity generalizes fairly easily because of its ``modular'' construction. In fact, all that we needed to know was that our inverse homomorphism $\alpha_{X} : K(\R^{2}\times X) \to K(X)$ satisfies the following two elementary properties:
\begin{description}
\item[Product Formula]
$\alpha_{X\times Y}(x\cdot y) = \alpha_{X}(x) \cdot y$
for all $x\in K(\R^{2}\times X)$ and $y\in K(Y)$.
\item[Normalization Axiom]
$\alpha_{\pt}(b) = 1$, where $b = \beta_{\pt}(1) \in K(\R^{2})$.
\end{description}
Hence, our proof would still hold if we were to give a different construction for the inverse homomorphism $\alpha_{X}: K(\R^{2}\times X ) \to K(X)$, provided that it still satisfies the above two properties. One other such construction is given in \cite{AB64}.

We now present an extension of our proof of Bott periodicity to the equivariant case. Suppose that we are given a compact Lie group $G$. Letting $V$ be a finite-dimensional complex $G$-module, we can construct the element $\lambda_{V}\in K_{G}(V)$ given by the external complex $\Lambda(V)$. Replacing $\R^{2}$ by $V$ in our discussion, we have a homomorphism
$$\varphi : K_{G}(X) \xrightarrow{ \lambda^{*}_{V}\tensorhat\cdot }
K_{G}(V \times X)$$
given by exterior multiplication by the dual bundle $\lambda^{*}_{V}$ of the Thom class $\lambda_{V}$. To show that this map is an isomorphism, we must then construct an inverse homomorphism $\alpha_{X}: K_{G}(V \times X) \to K_{G}(X)$ that satisfies the following two analogous axioms:
\begin{description}
\item[Product Formula]
$\alpha_{X\times Y}(x\cdot y) = \alpha_{X}(x) \cdot y$
for all $x\in K_{G}(V\times X)$ and $y\in K_{G}(Y)$.
\item[Normalization Axiom]
$\alpha_{\pt}(\lambda^{*}_{V}) = 1\in R(G)$, where $\lambda^{*}_{V} = \varphi_{\pt}(1) \in K_{G}(V)$.
\end{description}
Again, such a construction is possible using the index of a family of elliptic operators (see \cite[\S 4]{Ati68}), and so we obtain

\begin{theorem}[Equivariant Periodicity]
If $X$ is a compact $G$-space and $V$ is a finite dimensional complex $G$-module, then the homomorphism
$$\varphi : K_{G}(X) \to K_{G}(V \times X)$$
given by exterior multiplication by $\lambda_{V}^{*}\in K_{G}(V)$ is an isomorphism.
\end{theorem}

As a special case of this theorem, we can take $G = U(n)$ to be the unitary group and $V = \C^{n}$ to be the trivial representation. Then, any complex vector bundle over $X$ can be written in the form $E = Y \times_{U(n)}\C^{n}$, where  $Y$ is a free $U(n)$-space with $Y/U(n) = X$. Recalling from \S\ref{sec:4} that $K_{G}(X) \cong K(X/G)$ if $X$ is a free $G$-space, we see that
$$K_{U(n)}(Y) \cong K(X), \qquad
K_{U(n)}(Y \times_{U(n)} \C^{n}) \cong K(E).$$
Noting that the Thom class $\lambda_{E}\in K(E)$ corresponds to the element $\varphi(\lambda_{V}) \in K_{G}(V \times X)$, it follows immediately from the equivariant periodicity theorem that

\begin{corollary}[Thom Isomorphism Theorem]
If $E$ is a complex vector bundle over a compact space $X$, then multiplication by the Thom class $\lambda_{E}\in K(E)$ induces an isomorphism $K(X) \to K(E)$.
\end{corollary}

Furthermore, we can introduce the action of a second group $H$ without altering the discussion. If $X$ is an $H$-space, and $E$ is an $H$-bundle over $X$, then we can write $E$ in the form $E = Y \times_{U(n)} \C^{n}$, where $Y$ is a $H \times U(n)$-space with $Y/U(n) = X$. We then have
$$K_{H\times U(n)}(Y) \cong K_{H}(X), \qquad
K_{H\times U(n)}(Y \times_{U(n)}\C^{n}) \cong K_{H}(E),$$
and so we obtain the equivariant generalization.

\begin{corollary}[Equivariant Thom Isomorphism Theorem]
If $E$ is a complex $H$-vector bundle over a compact $H$-space $X$, then multiplication by the  Thom class $\lambda_{E}\in K_{H}(E)$ induces an isomorphism $K_{H}(X) \to K_{H}(E)$.
\end{corollary}

\section{The Atiyah-Hirzebruch Spectral Sequence}

Here we present the Atiyah-Hirzebruch spectral sequence for computing $K$-theory, as described in \cite[\S 2]{Ati-Hir}. We do not assume the reader is
familiar with the construction of a spectral sequence from a Cartan-Eilenberg system (see \cite{CE56}),
instead deriving all of the necessary homological algebra from the long exact sequences for pairs and triples.

Let $X$ be a finite (connected) CW-complex, and let $X^{p}$ denote its $p$-skeleton. Define
\begin{equation*}
K^{*}_{p}(X) = \Ker \bigl\{ K^{*}(X) \to K^{*}(X^{p-1})\bigr\}.
\end{equation*}
This is the Atiyah-Hirzebruch filtration for $K$-theory, with
$$K^{*}_{0}(X) = K^{*}(X) \,\supset\, K^{*}_{1}(X) = \tilde{K}^{*}(X)\, \supset \,\cdots \,\supset\, K^{*}_{\mathrm{max}} (X)= 0,$$
where $X^{-1} = \varnothing$ and $X^{\mathrm{max}-1} = X$.\footnote{Although it is not obvious from this definition, the Atiyah-Hirzebruch filtration is actually a homotopy type invariant, and it does not depend on the choice of CW-structure on $X$. Furthermore, this filtration is multiplicative, i.e., $K^{*}_{p}(X) \cdot K^{*}_{q}(X)  \subset K^{*}_{p+q}(X)$
for all $p$ and $q$. See \cite[\S 2]{Ati-Hir} for more information.} The Atiyah-Hirzebruch spectral sequence computes the associated graded components
$$\Gr_{p}K^{*}(X) = K^{*}_{p}(X)\,/\, K^{*}_{p+1}(X).$$ 
If we know these graded components, we can attempt to reconstruct $K^{*}(X)$ via a series of group extensions.

Consider now the following commutative diagram of exact sequences:
\begin{equation*}
\xymatrix@C=1em{
K^{*-1}(X^{p}) \ar[r] \ar[d] & K^{*}(X,X^{p}) \ar[r] \ar[d] & K^{*}(X) \ar[r] \ar@{=}[d] & K^{*}(X^{p}) \ar[d] \\
K^{*-1}(X^{p-1}) \ar[r] \ar[d] & K^{*}(X,X^{p-1}) \ar[r] \ar[d] & K^{*}(X) \ar[r] & K^{*}(X^{p-1}) \\
K^{*}(X^{p},X^{p-1}) \ar@{=}[r] \ar[d] & K^{*}(X^{p},X^{p-1}) \ar[d] \\
K^{*}(X^{p}) \ar[r] & K^{*+1}(X,X^{p})
}
\end{equation*}
where the top two rows are the long exact sequences for the pairs $(X,X^{p})$ and $(X,X^{p-1})$, and the left two columns are the long exact sequences for the pair $(X^{p},X^{p-1})$ and the triple $(X,X^{p},X^{p-1})$. From our definition 
of $K_{p}^{*}(X)$
and the exactness of the top two rows of the above sequence, we have in $K^{*}(X)$
\begin{align*}
K_{p}^{*}(X) &= \im  K^{*}(X,X^{p-1}) \cong K^{*}(X,X^{p-1})\,/\, \im K^{*-1}(X^{p-1}), \\
K_{p+1}^{*}(X) &= \im  K^{*}(X,X^{p}) \cong K^{*}(X,X^{p}) \,/\, \im K^{*-1}(X^{p}),
\end{align*}
and the associated graded component is thus
\begin{equation*}
\begin{split}
\Gr_{p}K^{*}(X) &\cong K^{*}(X,X^{p-1}) \, / \,
\Bigl( \im K^{*}(X,X^{p}) + \im K^{*-1}(X^{p-1}) \Bigr) \\
&\cong \Bigl( K^{*}(X,X^{p-1}) \,/\, \im K^{*}(X,X^{p}) \Bigr) \,/\, \im K^{*-1}(X^{p-1}) \\
&\cong \Ker \Bigl\{ K^{*}(X^{p},X^{p-1}) \to K^{*+1}(X,X^{p}) \Bigr\}\,/\, \im K^{*-1}(X^{p-1}).
\end{split}\end{equation*}
To compute this graded component, we construct a spectral sequence with first page,
\begin{equation*}
E^{p}_{1} := K^{*}(X^{p},X^{p-1})
\end{equation*}
and whose $k$-th page is given by:
\begin{equation*}
E^{p}_{k} = \frac{\Ker \Bigl\{ K^{*}(X^{p},X^{p-1}) \xrightarrow{\beta} K^{*+1}(X^{p-1+k},X^{p}) \Bigr\}}{ \im  \Bigl\{ K^{*-1}(X^{p-1},X^{p-k}) \xrightarrow{\beta'} K^{*}(X^{p},X^{p-1}) \Bigr\}}.
\end{equation*}
Note that for $E_{1}$, we start with
$$E^{p}_{1} = \Ker \Bigl\{ K^{*}(X^{p},X^{p-1}) \to K^{*+1}(X^{p},X^{p}) \Bigr\}\,/\, \im K^{*-1}(X^{p-1},X^{p-1}),
$$
which agrees with 
our definition of $E^{p}_{1}$.
On the other hand, for large enough $k$, we have $X^{p-1+k} = X$ and $X^{p-k} = \varnothing$,
and so the sequence stabilizes to give us our desired graded component $E_{\infty}^{p} = \Gr_{p}K^{*}(X)$. 
By examining the following commutative diagram, based on the bottom left corner of our previous diagram, 
$$
\xymatrix@C=1em{
K^{*-1}(X^{p-1},X^{p-k}) \ar[r] \ar[d] & K^{*}(X^{p-1+k},X^{p-1}) \ar[d] \\
K^{*}(X^{p},X^{p-1}) \ar@{=}[r] & K^{*}(X^{p},X^{p-1}) \ar[d] \\
& K^{*+1}(X^{p-1+k},X^{p})
}
$$
we see that in $K^{*}(X^{p},X^{p-1})$
\begin{equation*}\begin{split}
\im K^{*-1}(X^{p-1},X^{p-k})
&\subset \im K^{*-1}(X^{p-1+k},X^{p-1}) \\
&= \Ker \Bigl\{ K^{*}(X^{p},X^{p-1}) \to K^{*+1}(X^{p-1+k},X^{p}) \Bigr\},
\end{split}\end{equation*}
and thus $E^{p}_{k}$ 
is indeed well defined.

To construct the differentials for the spectral sequence, consider the following commutative diagram:
$$
\xymatrix@C=1em{
& K^{*-1}(X^{p-1},X^{p-k}) \ar[r] \ar[dr]_{\beta'} & K^{*-1}(X^{p-1},X^{p-k-1}) \ar[d]^{\alpha'} \\
&&K^{*}(X^{p},X^{p-1}) \ar@{-->}[dl]_{d_{k}} \ar[d]^{\alpha} \ar[dr]^{\beta} \\
K^{*}(X^{p-1+k},X^{p}) \ar[r]^-{\epsilon} & K^{*+1}(X^{p+k},X^{p-1+k})
\ar[r]^-{\delta} & K^{*+1}(X^{p+k},X^{p}) \ar[r]^-{\gamma} & K^{*+1}(X^{p-1+k},X^{p})
}
$$
where the top row comes from the long exact sequence for the triple $(X^{p-1},X^{p-k},X^{p-k-1})$ and the bottom row is the long exact sequence for the triple $(X^{p+k},X^{p-1+k},X^{p})$.
Our goal is to compute the groups $E_{k+1}^{p} = \Ker \alpha / \im \alpha'$ in terms of the groups $E_{k}^{p} = \Ker \beta / \im \beta'$.

Let $x \in K^{*}(X^{p},X^{p-1})$ be a representative of a class $[x]\in E_{k}^{p}$. Since $x\in\Ker\beta$ by the definition 
of $E_{k}^{p}$, it follows from our commutative diagram that $\alpha x \in \Ker\gamma = \im\delta$. We can therefore write $\alpha x = \delta y$ for some element $y\in K^{*+1}(X^{p+k},X^{p-1+k})$, where $y$ is determined up to $\Ker \delta = \im \epsilon$. This gives us a well defined class
$$[y] \in K^{*+1}\bigl(X^{p+k},X^{p-1+k}\bigr) \,/\,\im \Bigl\{ K^{*}\bigl(X^{p-1+k},X^{p}\bigr) \xrightarrow{\epsilon} K^{*+1}(X^{p+K},X^{p-1+k}) \Bigl\},$$
and we define
the differential of the spectral sequence to be
$$d_{k}: E^{p}_{k} \to E^{p+k}_{k}, \qquad d_{k} [x]:= [y].$$
To be thorough, one should further verify that
$$y \in \Ker \Bigl\{ K^{*+1}\bigl(X^{p+k},X^{p-1+k}\bigr) \to K^{*+2}\bigl(X^{p+2k-1},X^{p+k}\bigr)\Bigr\}$$
(by a similar diagram chase), that the class $d_{k}[x]\in E_{k}^{p+k,*+1}$ does not depend on the choice of $x$ (as $\im \beta' \subset \im \alpha' \subset \Ker \alpha$), and that the $d_{k}$-cohomology indeed gives
$$\Ker d_{k} / \im d_{k} = \Ker \alpha / \im \alpha' = E^{p}_{k+1},$$
the next page of the spectral sequence.

To insert the cohomology grading into the above discussion,
we use a bigrading on the spectral sequence, where the term $E^{p,q}_{k}$ carries cohomological degree $p+q$. The differential then yields
$d_{k} : E_{k}^{p,q} \to E_{k}^{p+k,q-k+1}$ for all $p$ and $q$.

The $E_{1}$ page of this spectral sequence has components $E_{1}^{p,q} = K^{p+q}(X^{p},X^{p-1})$.  However, if we work relative to the $(p-1)$-skeleton, then $X^{p}$ becomes a bouquet of $p$-spheres, with one for each $p$-cell. We therefore obtain
$$E_{1}^{p,q} = K^{p+q}\bigl(X^{p},X^{p-1}\bigr) \cong K^{p+q}(S^{p})^{\oplus \beta_{p}}
\cong K^{q}(\mathrm{pt})^{\oplus \beta_{p}}
= C^{p}\bigl(X;K^{q}(\mathrm{pt}) \bigr)$$
(where $\beta_{p}$ is the number of $p$-cells), giving the chain groups for the cellular cohomology of $X$ with coefficients in the $K$-theory of a point. In addition, the first differential
$$d_{1} : C^{p}\bigl(X;K^{q}(\mathrm{pt})\bigr) \to C^{p+1}\bigl(X;K^{q}(\mathrm{pt})\bigr)$$
is precisely the differential in the cellular chain complex (see \cite[\S III.7]{Ada74}), giving us
$$E_{2}^{p,q} = H^{p}\bigl(X;K^{q}(\mathrm{pt})\bigr)$$
as our starting point for the Atiyah-Hirzebruch spectral sequence.

The entire discussion above works not only for $K$-theory, but for any generalized cohomology theory.  In the specific case of complex $K$-theory, things simplify further.  We obtain
$$E_{2}^{p,q} = \begin{cases}
H^{p}(X;\mathbb{Z}) & \text{ for $q$ even, }\\
0 & \text{ for $q$ odd, }
\end{cases}$$
and so the entire spectral sequence can be computed in terms of the row with $q=0$.
All of the differentials $d_{k}$ with $k$ even vanish, as they interchange odd and even $q$.  On the other hand, for $k$ odd, the differentials are cohomology operations:
$$d_{k} : H^{p}(X;\mathbb{Z}) \to H^{p+k}(X;\mathbb{Z}).$$
So, to compute complex $K$-theory, we start with integral cohomology.  We then successively take the cohomology with respect to the operations $d_{k}$ for $k$ odd.  Finally, we recover $K^{0}(X)$ and $K^{1}(X)$ by combining what remains of all of the even and odd cohomology groups via a series of (often nontrivial) group extensions.  One consequence of this argument is that $K$-theory is never larger than integral cohomology.

\begin{example}
Consider the real projective space $X = \R P^{n}$. Its integral cohomology is: $$ H^{*}(\R P^{n};\Z) \cong \begin{cases}
\Z & \text{ in degree $0$,} \\
\Z/2\Z & \text{ in even degrees between $0$ and $n$,} \\
\Z & \text{ in degree $n$ if $n$ is odd,} \\
0 & \text{ otherwise.} \end{cases}
$$
The degree $0$ component $\Z$ must persist to $K$-theory, as that corresponds to the integral cohomology or $K$-theory of a point (and indeed this $\Z$ is not present in the reduced cohomology and $K$-theory), and thus all of the differentials $d_{k}$ vanish there. In addition, there are no non-zero homomorphisms $\Z/2\Z \to \Z$, and so the top degree component $\Z$ (for $n$ odd) always lies in the cokernel of the differentials $d_{k}$. Since the $\Z/2\Z$ components lie in only even degrees, it follows that all of the odd differentials $d_{k}$ vanish, and since the even differentials always vanish, we find that the Atiyah-Hirzebruch spectral sequence collapses at the $E_{2}$ page. It follows that $K^{*}(\R P^{n})$ has a component $\Z$ in degree 0, a component $\Z$ in degree 1 when $n$ is odd, and a torsion component in degree $0$ obtained by a series of $\Z/2\Z$-extensions of $\Z/2\Z$. However, this information does not completely determine the $K$-theory of $\R P^{n}$, as there are two possible forms for each of the $\Z/2\Z$-extensions:
\begin{align*}
 0 \to \Z/2\Z \to \Z/4\Z \to \Z/2\Z \to 0, \text{ or } \\
 0 \to \Z/2\Z \to \Z/2\Z \oplus \Z/2\Z \to \Z/2\Z \to 0.
\end{align*}
In fact, $K(\R P^{n})$ uses the former, giving us (see \cite[Proposition 2.7.7]{Ati67})
\begin{align*}
	K^{0}(\R P^{n} ) & \cong \Z \oplus \Z/2^{k}\Z, \\
	K^{1}(\R P^{n} ) &\cong \begin{cases}
		\Z & \text{ for $n$ odd, } \\
		0 & \text{ for $n$ even, } \\
	\end{cases}
\end{align*}
with $k = 2^{\lfloor n/2 \rfloor}$, as the $\Z/2\Z$ graded components successively combine via group extensions to give $\Z/2\Z$, $\Z/4\Z$, $\Z/8\Z$, etc.
\end{example}


\def\cprime{$'$}

\end{document}